\documentclass[12pt,leqno,draft]{article}
\usepackage{amsfonts}
\usepackage{amsmath, amsthm, amsfonts, amssymb, color}
\usepackage{mathrsfs}
\usepackage{color}
\newtheorem{thm}{Theorem}[section]

\newtheorem{rem}[thm]{Remark}
\theoremstyle{definition}

\newcommand{\scr}[1]{\mathscr #1}
\definecolor{wco}{rgb}{0.5,0.2,0.3}

\numberwithin{equation}{section} \theoremstyle{remark}

\newcommand{\ua}{\uparrow}

\title{{\bf   McKean-Vlasov SDEs with  Local Distributional Interactions: Well-Posedness and   Entropy-Cost Estimates}\footnote{Supported in part by the National Key R\&D Program of China (2022YFA1006000), NNSFC(12531007, 12301180, 12271398),  RGC(21301925), NSFC/RGC JRS N-CityU165/25 and  Research Centre for Nonlinear Analysis at Hong Kong PolyU.} }
\author{
{\bf   Xing Huang$^{a)}$,  Panpan Ren$^{b)}$, Feng-Yu Wang$^{a)}$  }\\
\footnotesize{ a) Center for Applied Mathematics and KL-AAGDM, Tianjin
University, Tianjin 300072, China}\\
\footnotesize{ b) Department of Mathematics, City University of  Hong Kong, Tat Chee Av., Hong Kong,  China }\\
\footnotesize{  xinghuang@tju.edu.cn, panparen@cityu.edu.hk, wangfy@tju.edu.cn}\\
}
\begin{document}
\allowdisplaybreaks
\def\R{\mathbb R}  \def\ff{\frac} \def\ss{\sqrt} \def\B{\mathbf
B} \def\W{\mathbb W}
\def\N{\mathbb N} \def\kk{\kappa} \def\m{{\bf m}}
\def\ee{\varepsilon}\def\ddd{D^*}
\def\dd{\delta} \def\DD{\Delta} \def\vv{\varepsilon} \def\rr{\rho}
\def\<{\langle} \def\>{\rangle} \def\GG{\Gamma} \def\gg{\gamma}
  \def\nn{\nabla} \def\pp{\partial} \def\E{\mathbb E}
\def\d{\text{\rm{d}}} \def\bb{\beta} \def\aa{\alpha} \def\D{\scr D}
  \def\si{\sigma} \def\ess{\text{\rm{ess}}}
\def\beg{\begin} \def\beq{\begin{equation}}  \def\F{\scr F}
\def\Ric{\text{\rm{Ric}}} \def\Hess{\text{\rm{Hess}}}
\def\e{\text{\rm{e}}} \def\ua{\underline a} \def\OO{\Omega}  \def\oo{\omega}
 \def\tt{\tilde} \def\Ric{\text{\rm{Ric}}}
\def\cut{\text{\rm{cut}}} \def\P{\mathbb P} \def\ifn{I_n(f^{\bigotimes n})}
\def\C{\scr C}      \def\aaa{\mathbf{r}}     \def\r{r}
\def\gap{\text{\rm{gap}}} \def\prr{\pi_{{\bf m},\varrho}}  \def\r{\mathbf r}
\def\Z{\mathbb Z} \def\vrr{\varrho} \def\ll{\lambda}
\def\L{\scr L}\def\Tt{\tt} \def\TT{\tt}\def\II{\mathbb I}
\def\i{{\rm in}}\def\Sect{{\rm Sect}}  \def\H{\mathbb H}
\def\M{\scr M}\def\Q{\mathbb Q} \def\texto{\text{o}} \def\LL{\Lambda}
\def\Rank{{\rm Rank}} \def\B{\scr B} \def\i{{\rm i}} \def\HR{\hat{\R}^d}
\def\to{\rightarrow}\def\l{\ell}\def\iint{\int}
\def\EE{\scr E}\def\Cut{{\rm Cut}}
\def\A{\scr A} \def\Lip{{\rm Lip}}
\def\BB{\scr B}\def\Ent{{\rm Ent}}\def\L{\scr L}
\def\R{\mathbb R}  \def\ff{\frac} \def\ss{\sqrt} \def\B{\mathbf
B}
\def\N{\mathbb N} \def\kk{\kappa} \def\m{{\bf m}}
\def\dd{\delta} \def\DD{\Delta} \def\vv{\varepsilon} \def\rr{\rho}
\def\<{\langle} \def\>{\rangle} \def\GG{\Gamma} \def\gg{\gamma}
  \def\nn{\nabla} \def\pp{\partial} \def\E{\mathbb E}
\def\d{\text{\rm{d}}} \def\bb{\beta} \def\aa{\alpha} \def\D{\scr D}
  \def\si{\sigma} \def\ess{\text{\rm{ess}}}
\def\beg{\begin} \def\beq{\begin{equation}}  \def\F{\scr F}
\def\Ric{\text{\rm{Ric}}} \def\Hess{\text{\rm{Hess}}}
\def\e{\text{\rm{e}}} \def\ua{\underline a} \def\OO{\Omega}  \def\oo{\omega}
 \def\tt{\tilde} \def\Ric{\text{\rm{Ric}}}
\def\cut{\text{\rm{cut}}} \def\P{\mathbb P} \def\ifn{I_n(f^{\bigotimes n})}
\def\C{\scr C}      \def\aaa{\mathbf{r}}     \def\r{r}
\def\gap{\text{\rm{gap}}} \def\prr{\pi_{{\bf m},\varrho}}  \def\r{\mathbf r}
\def\Z{\mathbb Z} \def\vrr{\varrho} \def\ll{\lambda}
\def\L{\scr L}\def\Tt{\tt} \def\TT{\tt}\def\II{\mathbb I}
\def\i{{\rm in}}\def\Sect{{\rm Sect}}  \def\H{\mathbb H}
\def\M{\scr M}\def\Q{\mathbb Q} \def\texto{\text{o}} \def\LL{\Lambda}
\def\Rank{{\rm Rank}} \def\B{\scr B} \def\i{{\rm i}} \def\HR{\hat{\R}^d}
\def\to{\rightarrow}\def\l{\ell}\def\BB{\mathbb B}
\def\8{\infty}\def\I{1}\def\U{\scr U} \def\n{{\mathbf n}}\def\v{V}
\maketitle

\begin{abstract}
We study   McKean-Vlasov SDEs with  interaction kernels in  $\tt W^{-\dd,k},$ the local negative Sobolev space on $\R^d$ with indexes $\dd \in  [0,\infty)$ and $k\in [1,\infty].$
We derive the  local well-posedness for  any singular indexes $(\dd,k)\in [0,\infty)\times [1,\infty],$ and prove the global well-posedness
 for any initial distributions  provided  $\dd+\ff d k<1$.   Moreover,    the relative entropy and the $\|\cdot\|_{\dd,k*}$-distance induced by  $ \tt W^{-\dd,k}$ are estimated for the time-marginal distributions of solutions by using the  Wasserstein distance of initial distributions, which describe the regularity of the solution in initial distribution. In particular, the main results apply to     Nemytskii-type SDEs  which    depend  on  higher order derivatives of the density functions, as well as McKean-Vlasov SDEs with interactions more singular   than Riesz kernels.

 \end{abstract} \noindent
 $\|\cdot\|_{\dd,k*}$-distance,   entropy-cost inequality.

\vskip 2cm

\tableofcontents

\section{Introduction}

Let $\scr P$ be the set of all probability measures on $\R^d$ equipped with the weak topology. Consider
the following  McKean-Vlasov SDE on $\R^d$:
\beq\label{E0} \d X_t= b_t(X_t, \L_{X_t})\d t+  \d W_t,\ \ t\ge 0,\end{equation}
where    $(W_t)_{t\ge 0}$ is a  $d$-dimensional Brownian motion on a probability base (i.e. complete filtered  probability space) $(\OO,\{\F_t\}_{t\ge 0},\F,\P)$, $\L_{X_t}$ is the distribution of $X_t$, $\tt{\scr P}$ is a measurable subspace of $\scr P$ to be determined by the singularity of $b_t(x,\mu)$ in $\mu$,  and
$$ b: [0,\infty)\times \R^d\times\tt{\scr P}\to \R^d$$
is measurable. Let $\B_b(\R^d)$ be the space of all bounded measurable functions on $\R^d$.

We are  interested in the case with singular interaction kernels where the drift includes
\beq\label{X} b_t(x,\mu)=(h_t*\mu)(x) := \int_{\R^d} h_t(x-y)\mu(\d y),\end{equation}
for $h_t$ belonging to a local negative Sobolev space, and the integral with respect to $\mu$ is understood as duality.

More precisely, let $(P_s^0:=\e^{s\DD})_{s\ge 0}$ be  the standard heat semigroup on $\R^d$, and recall that for any $\aa\in (0,\infty)$
\beq\label{FFD} (1-\DD)^{-\aa}:=\ff 1{\GG(\aa)} \int_0^\infty s^{\aa-1}\e^{-s}\,P_s^0\d s \end{equation}
is a bounded linear operator on $(\B_b(\R^d),\|\cdot\|_\infty)$. Let   $B(z,1):=\{x\in\R^d: |x-z|\le 1\}$ for $z\in\R^d$,   and
let    $\|\cdot\|_{L^k}$ be the $L^k$-norm with respect to the Lebesgue measure. Then
  for any  $\dd\in [0,\infty)$ and $k\in [1,\infty],$    the local negative Sobolev norm
$$\|f\|_{\tt W^{-\dd,k}}:= \sup_{z\in\R^d} \big\|1_{B(z,1)} (1-\DD)^{-\ff \dd 2} f\big\|_{L^k} $$
is well-defined on $\B_b(\R^d)$, and $\|\cdot\|_{\tt W^{-\dd,k}}\le c_{\dd,k} \|\cdot\|_\infty$ for some constant $c_{\dd,k}\in (0,\infty).$
 We define the local negative  Sobolev space $(\tt W^{-\dd,k},\|\cdot\|_{\tt W^{-\dd,k}})$ as the completion of $\B_b(\R^d)$ under the norm $\|\cdot\|_{\tt W^{-\dd,k}}$, which is
 a Banach space.

 When $\dd=0$, $(1-\DD)^{-\ff \dd 2}$ reduces to the identity operator so that
$$\big(\tt W^{-0,k},\|\cdot\|_{\tt W^{-0,k}}\big)= \big(\tt L^k, \|\cdot\|_{\tt L^k}\big),$$ where
$\tt L^k$ is the space of all functions $f\in L^k_{loc}(\R^d)$ with
$$\|f\|_{\tt L^k}:= \sup_{z\in\R^d} \big\|1_{B(z,1)}   f\big\|_{L^k}<\infty.$$  

Let $\tt W_*^{-\dd,k}$ be the dual space of $\tt W^{-\dd,k},$  which is a Banach space with    norm
 $$\|\mu-\nu\|_{\dd,k*}:=\sup_{f\in\B_b(\R^d), \|f\|_{ \tilde{W}^{-\dd,k}}\leq 1}   \big|\mu(f)-\nu(f)\big|,\ \ \ \mu,\nu\in \tt W^{-\dd,k}_*. $$
Then
$$\scr P_{\dd,k*}:= \scr P\cap \tt W_*^{-\dd,k}= \left\{\mu\in\scr P:\ \|\mu\|_{\dd,k*}:= \sup_{f\in \B_b(\R^d), \|f\|_{ \tilde{W}^{-\dd,k}}\leq 1}|\mu(f)|<\infty\right\}$$
is a   complete space under $\|\cdot\|_{\dd,k*},$ see Lemma \ref{L1} below.

For any
$$h_t=(h_t^i)_{1\le i\le d} \in \tt W^{-\dd,k}(\R^d;\R^d):= (\tt W^{-\dd,k})^d,$$
the drift  $b_t(x,\mu)$ in \eqref{X}  is well-defined   as the duality
$$b_t (x,\mu) := \ _{\tt W^{-\dd,k}}\big\<h_t(x-\cdot),\ \mu\big\>_{\tt W_*^{-\dd,k}}  =\Big(\ _{\tt W^{-\dd,k}}\big\<h_t^i(x-\cdot),\ \mu\big\>_{\tt W_*^{-\dd,k}}\Big)_{1\le i\le d} \in\R^d$$
for $  t\ge 0,\ x\in\R^d $ and $\mu\in \scr P_{\dd,k*}. $

 As a typical example of  interaction kernels,
the  Riesz kernel ${\bf K}:\R^d\to\R^d$  satisfies
$$|{\bf K}(z)|\le \ff c {|z|^\bb},\ \ 0\ne z\in\R^d$$ for some constants $c>0$ and $\bb\in (0,d)$, which belongs to $(\tt L^{k})^d= (\tt W^{-\dd,k})^d $ for $\dd=0$ and  $k\in [1,\ff d\bb)$. In this case,
 the well-posedness and regularity estimates have been established
in our recent paper \cite{HRW25} by establishing entropy-cost inequalities, which describe the regularity of time-marginal distributions of solutions with respect to initial distributions.
However, these regularity estimates remain  open when $h$ is merely distributional, although the well-posedness has  been intensively studied in \cite{CJM25,HRZ,IS} when $h$ belongs to a suitable Besov space, where the noise can be $\aa$-stable and/or degenerate. See also \cite{G,GP,HRZ,HZ} for the weak well-posedness  of  SDEs with
distributional  drifts without interaction (i.e. distribution independent), for which the entropy-cost inequality remains open.

In this paper,  we prove the well-posedness and establish regularity estimates for local distributional interaction kernels with arbitrary singular indexes
 $(\dd,k)\in [0,\infty)\times [1,\infty]$,
where the well-posedness result is also new in the literature, see Remark \ref{RM} below.

To solve \eqref{E0} for $b$ in \eqref{X} with $h_t\in \tt W^{-\dd,k}(\R^d;\R^d),$ we first present some  a-priori estimates on the time-marginal distribution of the solution, in
particular  we need to verify  that $\L_{X_t}\in  \scr P_{\dd, k*},$  which ensures   that the drift in \eqref{X} is well-defined for $\mu=\L_{X_t}.$
 To this end, we introduce the following estimate for the standard heat semigroup  $P_t^0$  (see  Lemma \ref{L0} below): there exists an increasing
 function $B: [0,\infty)\to (0,\infty)$ such that
 \beq\label{Hypin}\begin{split}&\|\nabla^iP_t^0\|_{\tt W^{-\dd,k}\to \tt W^{-\vv,p}}:=\sup_{\|f\|_{\tt W^{-\dd,k}}\le 1} \|P_t^0 f\|_{\tt W^{-\vv,p}}
  \le B_{\dd-\vv} t^{-\frac{i+\dd-\vv}{2}-\ff{d(p-k)}{2pk}},\\
  &\qquad \ \ t>0,\  i=0,1,\  \ \infty>\dd\ge \vv\ge 0,\ \ \infty\ge p\ge k\ge 1.
 \end{split}\end{equation}
When $i=0,$ this suggests  that the time-marginal distribution
$(\mu_t:=\L_{X_t})_{t\in  [0,T]}$ of solution  to \eqref{E0} satisfies
\begin{equation*} \rr_{\vv,p;\dd,k}^T(\mu):=\sup_{t\in [0,T]}t^{\frac{\dd-\vv}{2}+\ff{d(p-k)}{2pk}} \|\mu_t\|_{\dd,k*}\le C(T)\|\mu_0\|_{\vv,p*},\ \ \ T\in (0,\infty)\end{equation*}
for    some  increasing  $C: (0,\infty)\to  (0,\infty)$.  Thus, given initial value $X_0$ with
$\L_{X_0}\in \scr P_{\vv,p*}, $ it is reasonable to solve \eqref{E0} up to time $T$ with    time-marginal distributions belonging to the path space
  \beq\label{CPK}\C_{\vv,p;\dd,k}^{T}:=\Big\{\mu\in C^w([0,T];\scr P):\ \rr^{T}_{\vv,p;\dd,k}(\mu) <\infty\Big\},\end{equation}
where $C^w([0,T];\scr P) $ is the set of all weakly continuous maps from
$[0,T] $ to $  \scr P.$     This observation leads to the following  notion of  $\C_{\vv,p;\dd,k}$-solution of \eqref{E0}.

\beg{defn}[Maximal  $\C_{\vv,p;\dd,k}$-solution]  Let $1\le k\le  p\le \infty$, $0\le\vv\le\dd<\infty$, and
$\tt{\scr P}=\scr P_{\dd,k*}$.

\beg{enumerate}
\item[(1)]   We call $(X_t)_{t\in [0,\tau)}$ a  maximal strong  $\C_{\vv,p;\dd,k}$-solution of \eqref{E0} with life time $\tau$, if $\L_{X_0}\in \scr P_{\vv,p*}$, $\tau\in (0,\infty] $   such that
$$\limsup_{t\uparrow \tau} \|\L_{X_t}\|_{\dd,k*}=\infty \ \ \text{when }\ \tau<\infty,$$
$ (\L_{X_t})_{t\in [0,T]}  \in \C_{\vv,p;\dd,k}^T$ for any $T\in (0,\tau)$, and   $\P$-a.s.
$$X_t= X_0+\int_0^t b_s(X_s,\L_{X_s})\d s+ W_t,\ \ t\in [0,\tau).$$
When  $\tau=\infty$, we call $(X_t)_{t\ge 0}$  a  global strong $\C_{\vv,p;\dd,k}$-solution of \eqref{E0}.  For any $T\in (0,\tau)$, we call  $(X_t,W_t)_{t\in [0,T]}$   a   strong  $\C_{\vv,p;\dd,k}$-solution of \eqref{E0} up to time $T$.

\item[(2)]    A couple $(X_t,W_t)_{t\in [0,\tau)}$ is called a maximal weak $\C_{\vv,p;\dd,k}$-solution of \eqref{E0} with initial distribution $\gg\in \scr P_{\vv,p*}$,
 if there exists a probability base $(\OO,\{\F_t\}_{t\in [0,\tau)},\F,\P)$
such that $(W_t)_{t\in [0,\tau)}$ is a  $d$-dimensional Brownian motion, $\L_{X_0}=\gg$ and $(X_t)_{t\in [0,\tau)}$ is a maximal strong $\C_{\vv,p;\dd,k}$-solution of
\eqref{E0}. For any $ T\in (0,\tau)$, $(X_t,W_t)_{t\in [0,T]}$ is called a   weak $\C_{\vv,p;\dd,k}$-solution of \eqref{E0} up to time $T$.
\item[(3)]  If \eqref{E0} has a maximal weak $\C_{\vv,p;\dd,k}$-solution with initial distribution $\gg$,   and   any two maximal weak $\C_{\vv,p;\dd,k}$-solutions with initial distribution $\gg$ have common life time and distribution,
   then we say that \eqref{E0}    has a unique maximal weak $\C_{\vv,p;\dd,k}$-solution with initial distribution $\gg$. In this case, we denote the life time by $\tau(\gg)$, and set
$$P_t^*\gg:=\L_{X_t},\ \ t\in [0,\tau(\gg)).$$
 \end{enumerate}\end{defn}

In Section 2, we state the main results of the paper on the well-posedness and regularity estimates for $\C_{\vv,p;\dd,k}$-solutions of \eqref{E0}, see Theorem \ref{T0} and Theorem \ref{T02}. Section 3 contains necessary preparations, which will be  used in Sections 4 and 5 to prove  these two theorems respectively.

\section{Main results  }

\subsection{Well-posedness}

 To solve \eqref{E0},  we make the following assumptions where the drift is Lipschitz continuous in distribution under the $\|\cdot\|_{\dd,k*}$ distance.
 To cancel the singularity in small times caused by \eqref{Hypin}, we allow the drift vanishing at $t=0$   with rate $t^\kk$ for some $\kk\ge 0$.

\beg{enumerate} \item[{\bf (A)}] Let $1\le k \le \infty$  and $\dd, \kk\in [0,\infty)$. There exists increasing   $K:  (0,\infty)\to [1,\infty)$  such that
$$  |b_t(x,\nu) | \le K_t  t^\kk  \|\nu\|_{\dd,k*},\ \ \ \ |b_t(x,\mu)-b_t(x,\nu)|\le K_t t^\kk \|\mu- \nu\|_{\dd,k*}$$
hold for any $t\in (0,\infty),\ x\in \R^d,\ \mu,\nu\in \scr P_{\dd,k*}.$
\end{enumerate}

Under {\bf (A)}, let     $\vv\in [0,\dd]$ and $p\in [k,\infty]$ such that
\beq\label{TJ'} \eta:=\dd-\vv  +\ff{d(p-k)}{pk}<   1+2\kk.\end{equation}
 In this case,
 \beq\label{TH} \theta:=\ff 2 {1-(\eta-2\kk)^+}\in [2,\infty).\end{equation}

 \beg{thm}\label{T0}  Assume {\bf (A)}. Let  $\vv\in [0,\dd]$ and $p\in [k,\infty]$ satisfying $\eqref{TJ'}$.
 \beg{enumerate} \item[$(1)$]   For any $\F_0$-measurable initial value $X_0$ with $\gg:=\L_{X_0}\in  \scr P_{\vv,p*}$,   $\eqref{E0}$   has a unique maximal (weak and strong) $\C_{\vv,p;\dd,k}$-solution,
and there exists  increasing  $$C_\gg: [1,\infty)\times (0,\tau(\gg))\to (0,\infty)$$ such that
\beq\label{NES} \E\bigg[\sup_{s\in [0,t]} |X_s|^q\bigg|\F_0\bigg]\le C_\gg(q,t) (1+|X_0|^q),\ \ q\in [1,\infty),\ t\in (0,\tau(\gg)).\end{equation}
If $\vv=0,p=\infty$ and $\dd+\ff d k<1$, then $\tau(\gg)=\infty$ and   $C_\gg(q,t)=C(q,t)$ is independent of $\gg\in \scr P. $
\item[$(2)$]   Let $\theta$ be in $\eqref{TH}$.  For any   $n\in\mathbb N$, there exist constants $A_n \in (0,\infty)$  and $\ll_n\in [0,\infty)$ such that
for any  $\gg\in \scr P_{\vv,p*}$,
\beq\label{TT0} \tau(\gg)> \tau_n(\gg):=\beg{cases} n, &\text{if}\ \vv=0,p=\infty,\\
 \min\big\{n,\ \big(A_n \e^{A_n\|\gg\|_{\vv,p*}^{\theta}  }\big)^{-1}\big\},  &\text{otherwise},\end{cases} \end{equation}
\beq\label{EST}  \sup_{t\in (0,\tau_n(\gg)]} t^{\ff\eta 2}\e^{-\ll_n t}\|P_t^*\gg\|_{\dd,k*}\le 2B_{\dd-\vv}\|\gg\|_{\vv,p*},\end{equation}
where $\ll_n=0$ if $\tau_n(\gg)<n$.
In particular, if $\vv=0$ and $p=\infty$, then $\tau(\gg)=\infty$ for any $\gg\in \scr P$, and
\begin{equation*}\sup_{t\in (0,T],\gg\in \scr P} t^{\ff\eta 2} \|P_t^*\gg\|_{\dd,k*}<\infty,\ \ T\in (0,\infty).\end{equation*}
  \end{enumerate}
\end{thm}

 \begin{rem}\label{RM}  Let $b$ be in \eqref{X} such that
 $$\|h_t\|_{\tt W^{-\dd,k}}\le K_t t^\kk,\ \ \ t\ge 0.$$
Then assertions in Theorem $\ref{T0}$ hold for $\vv\in [0,\dd]$ and $p\in [k,\infty]$ satisfying \eqref{TJ'}.

 When $\kk=0$,  the condition  $\eqref{TJ'} $ coincides with    {\rm \cite[(1.33)]{HRZ}} for $\aa=2, q_b=\infty, p_0=\ff{p}{p-1},  \rho_0=k,  \beta_0=\vv$  and $ \beta_b=-\dd$. In this case,    {\rm \cite[Theorem 1.9]{HRZ}} ensures the well-posedness and density estimates for \eqref{E0} provided $\|h_t\|_{\BB_{k,\infty}^\dd}<\infty$, where $\BB_{k,\infty}^\dd$ is the Besov sapce. If moreover $p=\infty$ and $\vv=0$, \eqref{TJ'} becomes $\dd<1-\ff d k,$ so that  {\rm \cite[Theorem 1]{CJM25}} with $\alpha=2$
implies the global well-posedness of \eqref{E0} for
$\|h\|_{\BB_{k,\infty}^{-\dd}}<\infty.$
These results  do  not cover Theorem \ref{T0} since $\|h\|_{\BB_{k,\infty}^{-\dd}}<\infty $ may fail for $h\in \tt W^{-\dd,k}(\R^d;\R^d).$

For $\dd\in [0,\infty)$ and $k\in [1,\infty]$, let
$W^{-\dd,k}(\R^d;\R^d)$ be the closure of $C_0^\infty(\R^d;\R^d)$ with respect to the negative Sobolev norm
$\|f\|_{W^{-\dd,k}}:= \|(1-\DD)^{-\ff \dd 2} f\|_{L^k}$. The propagation of chaos is derived in {\rm \cite{JW}} for $h_t\in W^{-1,\infty}(\R^d;\R^d).$
However, the propagation of chaos remains open for   interactions
in $\tt W^{-\dd,k}(\R^d;\R^d) $ for $\dd>0$ and $k\ge 1$.

\end{rem}

\subsection{Regularity estimates }

Having  the maximal weak well-posedness for the $\C_{\vv,p;\dd,k}$-solution of \eqref{E0}, we aim to estimate   $\|P_t^*\gg-P_t^*\tt\gg\|_{\dd,k*}$ and the relative entropy $\Ent(P_t^*\gg|P_t^*\tt\gg),$
by  using the Wasserstein distance  $\W_q(\gg,\tt\gg)$ for some
$q\ge 1.$    Recall that for any $\gg,\tt \gg\in \scr P,$
$$\Ent(\gg|\tt\gg):= \beg{cases} \gg\big(\log\ff{\d\gg}{\d\tt\gg}\big),\ &\text{if}\ \ff{\d\gg}{\d\tt\gg}\ \text{exists},\\
\infty,\ &\text{otherwise},\end{cases}$$ and for any constant $q\in [1,\infty)$,
$$\W_q(\gg,\tt\gg):=\inf_{\pi\in \C(\gg,\tt\gg)}\bigg(\int_{\R^d\times\R^d} |x-y|^q\pi(\d x,\d y)\bigg)^{\ff 1 q},$$   where $\C(\gg,\tt\gg)$ is the set of all couplings for $\gg$ and $\tt\gg$.
Our estimates depend on
\beq\label{1ga} k_t(\gg):=  \|\gg\|_{\vv,p*}\lor \Big(\sup_{s\in (0,t]} s^{\frac{\dd-\vv}{2}+\ff{d(p-k)}{2pk}} \|P_s^*\gg\|_{\dd,k*}\Big),\ \ \ t\in (0,\tau(\gg)),\ \gg\in \scr P_{\vv,p*}.\end{equation}

Let $\theta$ be in \eqref{TH}. For any $\gg,\tt\gg\in \scr P_{\vv,p*}$  and  increasing function $\bb: (0,\infty)\to (0,\infty)$,  let
 \begin{equation*} K_{t,\bb}^{(\theta)}(\gg,\tt\gg):=
 \exp\Big[ \bb_t  \e^{\bb_t(tk_t(\gg)^{\theta}+ tk_t(\tt\gg)^{\theta})}\Big],\ \ \   t\in (0,\tau(\gg)\land\tau(\tt\gg)), \end{equation*}
\beq\label{ST}s_t(\theta',\gg):= t\land [k_t(\gg)^{-\theta'}],\ \ \ \theta'\in (\theta,\infty),\  t\in \big(0,\tau(\gg)\big).\end{equation}

 \beg{thm}\label{T02}   Assume {\bf (A)}. Let   $p\in [k,\infty]$ and $\vv\in [0, \dd\land \ff{d(p-1)}p]$ such that
 \beq\label{TJ} \eta:= \dd-\vv+\ff{d(p-k)}{pk}< 1\lor(\ff 12+\kk),\ \ \ \dd< 1\land \Big(2-\ff d k\Big)+(2\kk-\eta)^+.\end{equation}
 Then $1+(2\kk-\eta)^+-\eta>(\vv+\ff d p-\ff d k)^+$, and for any $q\in [1,\infty)$ satisfying
 \beq\label{QY} \ff{\vv+\ff d {p}}{1+(2\kk-\eta)^+-\eta}<q\le \ff{\vv+\ff d {p}}{(\vv+\ff d p-\ff d k)^+},\end{equation}
 where we set  $\ff{\vv+\ff d {p}}{(\vv+\ff d p-\ff d k)^+}=\infty$ if $\vv+\ff d p-\ff d k\le 0$, there exists  an  increasing function $\bb: [0,\infty)\to (0,\infty)$ such that the following assertions hold.
 \beg{enumerate} \item[$(1)$]    For any $\gg,\tt\gg\in \scr P_{\vv,p*}$ and  $t\in (0,\tau(\gg)\land \tau(\tt\gg)),$
  \beq\label{ES5} \beg{split}
 & \|P_t^*\gg-P_t^*\tt\gg\|_{\dd,k*}\\
 &\le (\|\gg\|_{\vv,p*}+\|\tt\gg\|_{\vv,p*})^{\ff{q-1}q} K_{t,\bb}^{(\theta)}(\gg,\tt\gg) t^{-[(\ff{1+\xi(q)}2)\vee(\frac{\delta}{2}+\frac{d}{2k})] }    \W_q(\gg,\tt\gg), \end{split}\end{equation}
 where, by  $\eqref{QY}$,
\beq\label{XQ} \xi(q):= \dd+\ff{dpq-k(d+\vv p)(q-1)}{pqk}\in \big[\eta,1+(2\kk-\eta)^+\big).\end{equation}
If $\vv=0, p=\infty$ and $\eta<1$,  then   for some increasing $\bb: (0,\infty)\to (0,\infty)$
\beq\label{ES5'}   \|P_t^*\gg-P_t^*\tt\gg\|_{\dd,k*} \le \bb_t    t^{-\ff {1+\dd} 2 -\ff d {2k}}\W_1(\gg,\tt\gg),\ \ t>0. \end{equation}
\item[$(2)$] For any $\theta'\in (\theta,\infty),$ $\gg,\tt\gg\in \scr P_{\vv,p*}$ and  $t\in (0,\tau(\gg)\land \tau(\tt\gg)),$
\beq\label{ES6} \beg{split}  \Ent(P_t^*\gg &|P_t^*\tt\gg)
 \le \bb_t(\|\gg\|_{\vv,p*}+\|\tt\gg\|_{\vv,p*})^{\ff{2(q-1)}q}\\
&\times   \bigg(\ff{\W_2(\gg,\tt\gg)^2}{s_t(\theta',\gg)}
+\ff  {K_{t,\bb}^{(\theta)}(\gg,\tt\gg)^2\W_q(\gg,\tt\gg)^2}{[s_t(\theta',\gg) \land s_t(\theta',\tt\gg)]^{([(1+\xi(q))\vee(\delta+\frac{d}{k})]-(2\kk+1))^+}}\bigg). \end{split}\end{equation}
 In particular, if   $\vv=0, p=\infty$ and $\eta<1$, then for some increasing $\bb: (0,\infty)\to (0,\infty)$
 \beq\label{ES7}  \Ent(P_t^*\gg|P_t^*\tt\gg)
\le   \ff {\bb_t}t  \W_2(\gg,\tt\gg)^2,\ \  t> 0,\ \gg,\tt\gg\in \scr P.  \end{equation}
\end{enumerate}
\end{thm}

 \beg{rem}To see that  $ \eqref{ES5}$ and $\eqref{ES6}$ characterize the regularity of the map $\gg\mapsto P_t^*\gg$, let
 $$P_tf(\gg):= \int_{\R^d} f\d (P_t^*\gg),\ \ \ f\in\B_b(\R^d).$$
 By Pinsker's inequality and \eqref{Var} below, there exists a constant $c\in (0,\infty)$ such that
 $$c\|\gg-\tt\gg\|_{var}\le \|\gg-\tt\gg\|_{\dd,k*}\land \ss{\Ent(\gg|\tt\gg)}.$$
 So,  each of $\eqref{ES5} $ and $\eqref{ES6}$ implies the local Lipschitz continuity of $\gg\mapsto P_tf(\gg)$ uniformly in $\|f\|_\infty\le 1$ with respect to $\W_q+\W_2$:
 $$\limsup_{\scr P_{\vv,p*} \ni \tt\gg\to\gg} \sup_{|f|\le 1} \ff{|P_tf(\gg)-P_tf(\tt\gg)|}{\W_2(\gg,\tt\gg)+\W_q(\gg,\tt\gg)}<\infty,\ \ \ \gg\in \scr P_{\vv,p*},\ t\in (0,\tau(\gg)).$$
 The estimates $\eqref{ES6}$-$\eqref{ES7}$ are called entropy-cost inequality or log-Harnack inequality. This type inequalities were first established in {\rm\cite{W10}} for elliptic diffusions on manifolds (possibly with boundary), see {\rm\cite{Wbook}} for the study of SPDEs, and see {\rm\cite{HW*,HRW,RW24}} for  the study of SDEs and McKean-Vlasov  SDEs.
 There are also many other papers concerning log-Harnack inequalities and applications, which we do not mention in details to save space.\end{rem}

We present the following  two examples to illustrate Theorem \ref{T0} and Theorem \ref{T02}, where the kernel $h$ is more singular than the Riesz kernel as considered in previous papers, see \cite{GLM,HRW25,L,RS,S} and references therein. In particular, our results apply to
 Nemytskii-type SDEs depending on higher order derivatives of the density.
 In the following example  the kernel $h$  is more singular than the Riesz kernel  ${\bf K}$ which satisfies   $|{\bf K}(z)|\le c |z|^{-\bb}$ for some $c \in (0,\infty)$ and $\bb\in (0,d)$.

\beg{exa}[Super Singular Interactions]      Let $b_t(x,\mu)=\int_{\R^d}h_t(x-y)\mu(\d y)$  for
$$h_t(z)=t^\kk h(z),\ \  h(z):= 1_{\{|z|>0\}}\ff {c z} {|z|^{d+1+\theta}},$$
where $\kk\in [0,\infty)$, $0\ne c\in\R,$   and $\theta\in  [0,1)$.  We have  {$h\in\tilde{W}^{-\dd,k}(\R^d,\R^d)$} for any   $k\in [1,\infty)$ and
\beq\label{**} \dd \in  {\Big(2+\ff{d(k-1)}{k},\ \infty\Big)},\end{equation}
so that   assertions in Theorem \ref{T0} hold for any $\vv\in [0,\dd]$ and $p\in [k,\infty]$ satisfying \eqref{TJ'},  and
assertions in Theorem   \ref{T02} hold under \eqref{TJ} and $\vv\le \ff{d(p-1)}{p}.$
  \end{exa}

\beg{proof}   By the definition of $\tt W^{-\dd,k}(\R^d,\R^d)=(\tt W^{-\dd,k})^d$,
it suffices to prove that the family
$$\{h_n:= 1_{\{|\cdot|\ge n^{-1}\}} h\}_{n\ge 1}\subset {\tilde{W}^{-\dd,k}(\R^d,\R^d)}$$ is a Cauchy sequence under $\|\cdot\|_{\tt W^{-\dd,k}}$, i.e.
\beq\label{CH} \lim_{n\to\infty} \sup_{m\ge n} \|h_n-h_m\|_{\tt W^{-\dd,k}}=0,\end{equation}
so that $h$ is well-defined in $ {\tilde{W}^{-\dd,k}(\R^d,\R^d)}$ as the limit of $h_n$   when $n\to\infty$.

To this end, we use \eqref{FFD} and the formula
$$P_t^0(h_n-h_m)(x)= (4\pi t)^{-\ff d 2} \int_{\{m^{-1}\le |z|<n^{-1}\} } h(z) \e^{-\ff{|z-x|^2}{4t}}\d z,\ \ \  t>0,\ x\in\R^d,\ n\le m.$$
By the integral transform $z\mapsto -z$ and $h(-z)=-h(z)$, we obtain
$$P_t^0(h_n-h_m)(x)= - (4\pi t)^{-\ff d 2} \int_{\{m^{-1}\le |z|<n^{-1}\} } h(z) \e^{-\ff{|z+x|^2}{4t}}\d z,\ \ \  t>0,\ x\in\R^d.$$
So, there exists a constant $c_1\in (0,\infty)$ such that
\beg{align*}&   \big|P_t^0(h_n-h_m)(x)\big| =\ff 1 2 \bigg|(4\pi t)^{-\ff d 2} \int_{\{m^{-1}\le |z|<n^{-1}\} } h(z) \Big(\e^{-\ff{|z-x|^2}{4t}}-\e^{-\ff{|z+x|^2}{4t}}\Big)\d z\bigg|\\
&\le \ff 1 2 (4\pi t)^{-\ff d 2}  \int_{\{m^{-1}\le |z|<n^{-1}\} }  \ff{| h(z)|\cdot |z|\cdot |x|} t  \Big(\e^{-\ff{|z-x|^2}{4t}}+ \e^{-\ff{|z+x|^2}{4t}}\Big) \d z\\
&= (4\pi t)^{-\ff d 2}  \int_{\{m^{-1}\le |z|<n^{-1}\} }  \ff{| h(z)|\cdot |z|\cdot |x|} t   \e^{-\ff{|z-x|^2}{4t}} \d z\\
&\le c_1   t^{-\ff d 2-1} \int_{B(0,n^{-1})}  \ff{  |x|} {|z|^{d+\theta-1}}  \e^{-\ff{|z-x|^2}{4t}} \d z.  \end{align*}
Noting that $|x|\le |z-x|+1$ for $|z|\le n^{-1}$,  by H\"older's inequality,  we find a constant $c_2\in (0,\infty)$ such that
\beg{align*} & t^{k+\ff{d k}2} \big\|   P_t^0(h_n-h_m)\big\|_{L^k}^k \le  \int_{\R^d}   \bigg(c_1  \int_{B(0,n^{-1})}   \ff{  |x|} {  |z|^{d+\theta-1}}  \e^{-\ff{|z-x|^2}{4t}} \d z\bigg)^k\d x\\
& \le  c_1^k \bigg(\int_{B(0,n^{-1})}  \ff{ 1} { |z|^{d+\theta-1}}  \d z \int_{\R^d}         (|z-x|+1)^k \e^{-\ff{k|z-x|^2}{4t}} \d x \bigg)\\
&\qquad\times  \bigg(  \int_{B(0,n^{-1})}       \ff{ 1} { |z|^{d+\theta-1}}  \d z\bigg)^{k-1} \\
&\le c_2^k n^{-k(1-\theta)}  \big(1+ t^{\ff k 2}\big) t^{\ff {d } 2}.\end{align*}
 Combining this with the formula \eqref{FFD}
  and noting that \eqref{**} implies
$$\ff \dd 2-2 -\ff{d(k-1)}{2k}  >-1,$$
we find a constant  $c_3 \in (0,\infty)$ such that
\beg{align*} &\sup_{m\ge n} \big\|h_n-h_m\big\|_{\tt W^{-\dd,k}}\le c(\dd)\int_0^\infty t^{\ff \dd 2-1} \e^{-t}\big\|   P_t^0(h_n-h_m)\big\|_{\tt L^k} \d t\\
&\le c_2 n^{-(1- \theta)}  \int_0^\infty t^{\ff \dd 2-2}   \big(1+ t^{\ff k 2}\big)^{\ff 1 k}  t^{-\ff{d(k-1)}{2k}   } \e^{-t}\d t\le c_3 n^{-(1-\theta)},\ \ n\ge 1. \end{align*}
This implies \eqref{CH} since $\theta\in (0,1)$.
  \end{proof}

Next, we consider SDEs whose coefficients depend on higher  order derivatives of the density function,   which include  the Burgers/Navier-Stokes/$p$-Laplacian equations
 as typical examples where the first and second order derivatives of density are involved, see for instance \cite{BRR,Szn}.

Let ${\bf \dd}_0$ be the Dirac function. Then for any absolutely continuous probability measure $\mu$ on $\R^d$, its density function can be formulated as
\beq\label{BAW0} \rr_\mu(x) = ({\bf \dd}_0*\mu)(x),\ \ \ x\in\R^d.\end{equation}
For any $i\in \mathbb N$, it is classical that  each component of $\nn^i {\bf \dd}_0$ belongs to
 $$ H^{-\dd}\subset \tt W^{-\dd,2}\ \ \text{if}\ \  \dd>\ff d 2+i,$$  so that for any $\mu\in \scr P_{\dd, 2*}$,
\beq\label{BAW1} \nn^i \rr_\mu(x):= (-1)^{i} _{\tt W^{-\dd, 2}}\big\<\nn^i {\bf \dd}_0(x-\cdot),\ \mu\big\>_{\tt W^{-\dd,2}_*} \end{equation}
 is well-defined in $\scr T_i:= \otimes^i \R^d,$   the space of $i$-tensors over  $\R^d$. When $\rr_\mu$ is regular enough, $\nn^i \rr_\mu$ defined in \eqref{BAW1} coincides with the corresponding classical derivatives.

  Given $n\in \mathbb N$, let  $\scr T_0:= \R$ and
$$\H_n:=\prod_{i=0}^{n-1} \scr T_i,$$ which
 is a finite-dimensional Hilbert space with induced norm $\|\cdot\|_{\H_n}.$
For any $\mu\in \scr P_{\dd, 2*}$ with $\dd> n-1+\ff d 2$ such that $\nn^i \rr_\mu$ exists for $0\le i\le n-1$,    denote
\beq\label{BAW2}  \rr_\mu^{\<n}(x):= (\nn^{i} \rr_\mu (x))_{0\le i\le n-1}\in \H_n,\ \ \ x\in \R^d.\end{equation}
In particular, when $n=1$ we have $\rr_\mu^{\<1}=\rr_\mu$. Moreover, let $\ell_{\xi}$ denote the distribution density  function for an absolutely continuous random variable $\xi$ on $\R^d$.

Now,   we consider the following SDE on $\R^d$ for a fixed time $T>0$:
\beq\label{E*}   \d X_t=    \d W_t
    +b_t\big(X_t, \ell_{X_t}^{\<n}(X_t)\big)\d t,\ \ t\in [0,T],\end{equation}
where
$$  b: [0,T]\times \R^d\times \H_n\to \R^d$$
is measurable.

\beg{exa}[Density-Derivative Dependent SDE]   If there exist $\kk\in [0,\infty)$ and increasing $K: (0,\infty)\to [1,\infty) $ such that
  \beg{align*}&  |b_t(x,h)-b_t(x,\tt h)|\le K_t t^\kk \|h-\tt h\|_{\H_n},\\
&|b_t(x,h)|\le K_t t^\kk (1+\|h\|_{\H_n}),\ \ \ t\in [0,\infty),\ h,\tt h\in \H_n,\ x\in\R^d.\end{align*}
  Then by \eqref{BAW0}-\eqref{BAW2}, {\bf (A)} holds for $k=2$ and any $\dd>\ff d2+n-1.$
So,   for the density-derivative  dependent SDE \eqref{E*},   assertions in Theorem \ref{T0} hold for any $\vv\in [0,\dd]$ and $p\in [2,\infty]$ satisfying \eqref{TJ'},  and
assertions in Theorem   \ref{T02} hold under \eqref{TJ} and $\vv\le \ff{d(p-1)}{p}.$
In particular, if $1+2\kk> d+n-1$, then we may take $\vv=0$ and $p=\infty$ such that  \eqref{E*} has   a unique (weak and strong) global $\C_{0,\infty;\dd,2}$-solution for
  {$\dd\in(\frac{d}{2}+n-1,1+2\kappa-\frac{d}{2})$} and any initial distribution $\mu\in \scr P$,
 and  when $\ff 1 2+\kk>d+n-1$ there exists increasing $\bb: (0,\infty)\to (0,\infty)$ such that
$$\Ent(P_t^*\mu|P_t^*\nu)\le \ff{\bb_t}t \W_2(\mu,\nu)^2,\ \ \ t>0,\ \mu,\nu\in  \scr P.$$
 \end{exa}

\section{Some preparations}

We will frequently use the following  simple  inequality: for any $\aa_1,\aa_2\in [0,1)$ and $\aa\in [0,1-\aa_2],$ there exists $c(\aa,\aa_1,\aa_2)\in (0,\infty)$ such that
\beq\label{LN0} \int_0^t s^{-\aa_1}(t-s)^{-\aa_2}\e^{-\ll (t-s)}\d s\le c(\aa,\aa_1,\aa_2) t^{1-\aa-\aa_1-\aa_2}\ll^{-\aa}, \ \ \ t,\ll>0.  \end{equation}
Indeed, by  the FKG inequality, when $ \aa_1,\aa_2\in [0,1)$, we have
\beg{align*} \int_0^t  {s^{-\aa_1}}(t-s)^{-\aa_2}\e^{-\ll (t-s)}\d s&\le \bigg(\ff 1 t\int_0^t  {s^{-\aa_1}}\d s\bigg)\int_0^t (t-s)^{-\aa_2}\e^{-\ll (t-s)}\d s \\
&= \ff 1 {1-\aa_1} t^{-\aa_1}\int_0^t s^{-\aa_2}\e^{-\ll s}\d s,\ \  t,\ll>0.\end{align*}
Then for $\aa\in [0,1-\aa_2)$ the inequality \eqref{LN0} follows from H\"older's inequality
\beg{align*} &\int_0^t s^{-\aa_2}\e^{-\ll s}\d s   \le   \bigg(\int_0^t s^{-\ff{\aa_2}{1-\aa}}\d s \bigg)^{1-\aa} \bigg(\int_0^t \e^{-{\frac{\ll}{\aa}} s}\d s\bigg)^{\aa}\\
&\le {\aa^\aa }\Big(\ff{1-\aa}{1-\aa-\aa_2}\Big)^{1-\aa} t^{{1-\aa-\aa_2}}\ll^{-\aa},\ \
  \ t,\ll>0,\end{align*}
and when $\aa=1-\aa_2$, \eqref{LN0} is implied by
$$  \int_0^\infty s^{-\aa_2}\e^{-\ll s}\d s \le \int_0^{\ll^{-1}}s^{-\aa_2}\d s+\ll^{\aa_2} \int_{\ll^{-1}}^\infty \e^{-\ll s}\d s\le \ff{2-\aa_2}{1-\aa_2} \ll^{\aa_2-1},\ \ \ll>0.$$

\beg{lem}\label{L0} Let $P_t^0$ be the   heat semigroup generated by $\DD$ on $\R^d$. Then for any $i_0\in\mathbb N$,  there exists  increasing $B: [0,\infty)\to (0,\infty)$ such that
\beg{align*}&\|\nn^i P_t^0\|_{\tt W^{-\dd,k}\to\tt W^{-\vv,p}}\le B_{\dd-\vv} t^{-\ff {i+\dd-\vv}2 -\ff{d(p-k)}{2pk}},\\
&\  \ t>0,\ 0\le \vv\le \dd<\infty,\ 1\le k\le p\le \infty,\ 0\le i\le i_0.\end{align*}
\end{lem}

\beg{proof} Noting that
$$P_t^0 f(x)= \big(4\pi t\big)^{-\ff d 2} \int_{\R^d} \e^{-\ff{|x-y|^2}{4t}} f(y)\d y,$$
and there exists a constant  $L\in (0,\infty)$ such that
$$\Big|\nn^i \e^{-\ff{|\cdot|^2}{4t}}\Big|(x)\le  L  t^{-\ff i 2} \e^{-\ff{|x|^2}{8t}},\ \ \ t>0,\ 0\le i\le i_0,\ x\in\R^d,$$
we  find an increasing function $C: \Z_+\to (0,\infty)$ such that
$$|(1-\DD)^n \nn^i P_t^0f|\le C_{n} t^{-\ff i 2 -n} P_{2t}^0 |f|,\ \ t>0,\ 0\le i\le i_0,n\in \Z_+, f\in \B_b(\R^d).$$
On the other hand, by \cite[Lemma 5.3]{HRW25}, there exists a constant $c_0\in (0,\infty)$ such that
$$\|P_{2t}^0\|_{\tt L^k\to\tt L^p}\le c_0 t^{-\ff{d(p-k)}{2pk}},\ \ t>0,\ 1\le k\le p\le \infty.$$
So,
\beq\label{PL*} \|(1-\DD)^n \nn^i P_t^0\|_{\tt L^k\to\tt L^p}\le C_{n} c_0 t^{-\ff i 2-n-\ff{d(p-k)}{2pk}},\ \ t>0,\ 0\le i\le i_0,n\in \Z_+.\end{equation}
 Let $n\in \mathbb N$ such that
$$\theta_0:=n-\ff{\dd-\vv}2\in [1,2).$$
By \eqref{FFD} and \eqref{PL*},
 we find an increasing  function   {$\bar{B}: \Z_+\to  (0,\infty)$} such that
\beg{align*} &\|\nn^i P_t^0\|_{\tt W^{-\dd,k}\to\tt W^{-\vv,p}}=\|(1-\DD)^{-\ff\vv 2}\nabla^iP_t^0 (1-\DD)^{\ff\dd 2}\|_{\tt L^k\to\tt L^p} \\
&= \|(1-\DD)^n\nn^i(1-\DD)^{-\theta_0}P_t^0\|_{\tt L^k\to\tt L^p}\\
&=\ff 1 {\GG(\theta_0)} \bigg\|\int_0^\infty s^{\theta_0-1}\e^{-s} (1-\DD)^n\nn^i P_{t+s}^0\d s\bigg\|_{\tt L^k\to\tt L^p}\\
&\le \ff {C_{n}c_0} {\GG(\theta_0)} \int_0^\infty s^{\theta_0-1}\e^{-s} (t+s)^{-\ff i 2-n-\ff{d(p-k)}{2pk}}\d s \\
&\le {\bar{ B}_{n}} t^{-\ff{i+\dd-\vv}2 - \ff{d(p-k)}{2pk}},\ \ t>0,\ 0\le i\le i_0,\ \delta-\vv\in (2(n-2), 2(n-1)],\ n\in\mathbb N. \end{align*}
This implies the desired estimate for some increasing function $B: [0,\infty)\to (0,\infty).$
\end{proof}

\beg{lem}\label{L1} Let $1\le k\le p\le\infty$,  $\varepsilon\le\dd <\infty$,    $\ll\in [0,\infty)$ and $T\in (0,\infty)$.
\beg{enumerate} \item[$(1)$]   The metric space $(\scr P_{\dd,k*},\ \|\cdot\|_{\dd,k*})$  is   complete, and the Borel $\si$-field coincides with that induced by the weak topology.
\item[$(2)$] For any $\gg\in \scr P_{\vv,p*},$
 the space $(\C_{\vv,p;\dd,k}^{\gg,T},\rr^{\ll,T}_{\vv,p;\dd,k})$ is complete, where for $\C_{\vv,p;\dd,k}^{T}$ in $\eqref{CPK}$,
  \beg{align*}& \C_{\vv,p;\dd,k}^{\gg,T}:=\big\{\mu \in \C_{\vv,p;\dd,k}^{T}:\ \mu_0=\gg\big\},\\
  &\rr^{\ll,T}_{\vv,p;\dd,k}(\mu,\nu):=\sup_{t\in (0,T]} \e^{-\ll t} t^{\frac{\dd-\vv}{2}+\ff{d(p-k)}{2pk}} \|\mu_t-\nu_t\|_{\dd,k*}.\end{align*}
\end{enumerate}  \end{lem}

\beg{proof}  (1) To prove the completeness of $(\scr P_{\dd,k*},\ \|\cdot\|_{\dd,k*})$, let  $\{\mu_n\}_{n\ge 1} $ be a Cauchy sequence in $(\scr P_{\dd,k*},\ \|\cdot\|_{\dd,k*})$. Since $(\scr P_{\dd,k*},\ \|\cdot\|_{\dd,k*}) $ is included by the dual space  $(\tt W^{-\dd,k}_*,\ \|\cdot\|_{\dd,k*})$ of the Banach space $(\tt W^{-\dd,k},\|\cdot\|_{\tt W^{-\dd,k}})$, there exists a unique $\mu\in \tt W^{-\dd,k}_*$ such that
$$\lim_{n\to\infty} \|\mu_n-\mu\|_{\dd,k*}=\sup_{f\in\B_b(\R^d), \|f\|_{\tt W^{-\dd,k}}\le 1} |\mu_n(f)-\mu(f)|=0.$$
Since $\scr P_{\dd,k*}= \tt W_*^{-\dd,k}\cap \scr P,$
it remains to show that $\mu\in \scr P$.  When $\dd=0$, we have
$$\|f\|_{\tt W^{-0,k}}=\|f\|_{\tt L^k}\le \oo(d)^{\ff 1 k} \|f\|_\infty,$$
where $\oo(d)$ is the volume of the unit ball $B(0,1)$. If $\dd>0$, by \eqref{FFD} we have
\beg{align*} \|f\|_{\tt W^{-\dd,k}}&=\bigg\|\ff1{\GG(\dd/2) } \int_0^\infty t^{\ff\dd 2-1} \e^{-t} P_t^0f\d t\bigg\|_{\tt L^k}\\
&\le \ff{\oo(d)^{\ff 1 k}}{\GG(\dd/2) }   \bigg\|  \int_0^\infty t^{\ff\dd 2-1} \e^{-t} P_t^0f\d t\bigg\|_{\tt L^\infty}\\
&\le \bigg(\ff{\oo(d)^{\ff 1 k}  }{\GG(\dd/2) } \int_0^\infty t^{\ff\dd 2-1} \e^{-t} \d t \bigg) \|f\|_{\infty}.\end{align*}
In any case, we find a constant $c\in (0,\infty)$ such that $\|\cdot\|_{\tt W^{-\dd,k}}\le c \|\cdot\|_\infty$, hence
\beq\label{Var} \|\mu_n-\mu\|_{var} \le  c \|\mu_n-\mu\|_{\dd, k*},\end{equation} so that
$$\lim_{n\to\infty} \|\mu_n-\mu\|_{var} \le  c \lim_{n\to\infty}\|\mu_n-\mu\|_{\dd, k*}  =0.$$
This together with  $\{\mu_n\}_{n\ge 1}\subset \scr P$ implies $\mu\in \scr P$.

 Next, since $C_b(\R^d)$ is dense in $\tt W^{-\dd,k}$,  {then for any $f\in    \tt W^{-\dd,k}$ and $\mu\in\scr P_{\delta,k*}$, there exists $\{f_n\}_{n\geq 1}\subset C_b(\R^d)$ satisfying  $\lim_{n\to\infty}\|f_n-f\|_{\tt W^{-\dd,k}}= 0$, which implies
 $$\lim_{n\to\infty}|\mu(f_n)-\mu(f)|\leq \|\mu\|_{\delta,k*}\lim_{n\to\infty}\|f_n-f\|_{\tt W^{-\dd,k}}= 0.$$}
 Noting that the Borel $\sigma$-field in $\scr P_{\dd,k*}$ is induced by
 $$\big\{\mu\mapsto \mu(f):\ f\in    \tt W^{-\dd,k}\big\},$$
 hence it coincides with the $\sigma$-field induced by the weak topology.

(2) It suffices to prove for $\ll=0.$   Let $\{\mu^{(n)}\}_{n\ge 1}\subset \C_{\vv,p;\dd,k}^{\gg,T}$ be
a Cauchy sequence with respect to $\rr^{T}_{\vv,p;\dd,k}$. Then for any $t\in (0,T]$, $\{\mu_t^{(n)}\}_{n\ge 1}$ is a Cauchy sequence
in $\scr P_{\dd,k*}$, so that by (1), there exists a unique $\mu_t\in \scr P_{\dd,k*}$ such that
\beq\label{IO}  \beg{split}&\lim_{n\to\infty}  \rr^{T}_{\vv,p;\dd,k}(\mu^{(n)}, \mu )
   =  \lim_{n\to\infty} \sup_{t\in [0,T]} \lim_{l\to\infty}  t^{\ff{\dd-\vv}2 +\ff{d(p-k)}{2pk} }\|\mu_t^{(n)} -\mu_t^{(l)} \|_{\dd,k*} \\
  & \le   \lim_{n,l\to\infty} \rr^{T}_{\vv,p;\dd,k}(\mu^{(n)}, \mu^{(l)})=0.  \end{split}\end{equation}
It remains to show the weak continuity of $(0,T]\ni t \mapsto \mu_t$, which together with \eqref{IO} implies   $\mu\in \C_{\vv,p;\dd,k}^{\gg,T}.$
 For any $f\in C_b(\R^d), s\in (0,T]$ and $\vv'>0$, by \eqref{Var}  and applying \eqref{IO}, we find large enough $n\ge 1$ such that
\beg{align*}&|\mu_t^{(n)}(f)-\mu_t(f)| \le c\|f\|_\infty \|\mu_t^{(n)}- \mu_t\|_{\dd,k*} \\
&\le c\|f\|_\infty   t^{-\ff{\dd-\vv}2-\ff{d(p-k)}{2pk}} \rr^{T}_{\vv,p;\dd,k}(\mu^{(n)}, \mu ) \le \vv',\ \ t\in [s/2,T].\end{align*}
Combining this with the weak continuity of $(0,T]\ni t\mapsto \mu_t^{(n)},$  we derive
$$\limsup_{t\to s} |\mu_t (f)-\mu_s(f)| \le \limsup_{t\to s} \big\{|\mu_t^{(n)} (f)-\mu_s^{(n)} (f)|+2\vv'\big\}=2\vv'.$$
Since $\vv'>0$ is arbitrary, this implies that weak continuity of $(0,T]\ni t \mapsto \mu_t.$
 \end{proof}

For any $T\in (0,\infty)$ and $\mu\in\C_{\vv,p;\dd,k}^{T}$, consider the SDE
\beq\label{BX} \d   X_{s,t}^{\mu,x}= b_t(  X_{s,t}^{\mu,x},\mu_t)\d t+ \d W_t,\ \ t\in [s,T],\   X_{s,s}^{\mu,x}=x.\end{equation}

 \beg{lem}\label{L2} Assume {\bf(A)} and let $\vv\in [0,\dd]$ and $p\in [k,\infty] $ satisfying $\eqref{TJ'}$.
  Then for any  $T\in (0,\infty)$ and  $\mu\in\C_{\vv,p;\dd,k}^{T}$, the SDE $\eqref{BX}$ is (weakly and  strongly) well-posed.
 \end{lem}

 \beg{proof}   By {\bf(A)} and $\mu\in \C_{\vv,p;\dd,k}^{T}$, there exists a constant $c\in (0,\infty)$ such that
 $b_t^{\mu}(x):=b_t(x,\mu_t)$ satisfies
 $$|b_t^{\mu}(x)|\le c t^{\kk -\ff\eta 2},\ \ t\in (0,T].$$
 By this and \eqref{TJ'}, we find $q'>2$ such that
\beq\label{ZV} \|b^{\mu}\|_{\tt L_{q'}^{\infty}(T)}=\sup_{z\in\R^d}\left(\int_{0}^T\|b_t^\mu 1_{B(z,1)}\|_{\infty}^{q'}\right)^{\frac{1}{q'}}<\infty.\end{equation}
Then the desired assertion follows from \cite[Proposition 5.1]{HRW25}.
 \end{proof}

For $\mu\in \C_{\vv,p;\dd,k}^{T}$, under \eqref{TJ'} we denote
\beq\label{SM}  P_{s,t}^\mu f(x):= \E[ f(  X_{s,t}^{\mu,x})],\ \ 0\le s\le t\le T,\ f\in \B_b(\R^d), x\in\R^d.\end{equation}

The next lemma provides the estimate \eqref{Hypin} for $ {P}_{s,t}^\mu$ replacing $P_t^0$, which is crucial in the proof of the main results.

\beg{lem}\label{L3'} Assume {\bf(A)} and let $\vv\in [0,\dd]$ and $p\in [k,\infty] $ satisfying $\eqref{TJ'}$.   Let  $\theta $ be in  $\eqref{TH},$ $p_1\in [1,\infty],   p_2 \in [p_1,\infty]$,  $\infty>\vv_1\ge \vv_2\ge 0,$   and $i=0,1.$
\beg{enumerate} \item[$(1)$]  If $\xi:=\vv_1-\vv_2 +\ff{d(p_2-p_1)}{p_1p_2}<1\lor\big(2-i+2\kk-\eta\big)$,
 then
  there exists an increasing function    $\beta: [0,\infty)\to (0,\infty)$  such that    for any  $ t\in (0,\infty) $   and   $  \mu\in \C_{\vv,p; \dd,k}^{t},$
\beq\label{ES1'} \|\nn^i  P_{t}^\mu\|_{\tt W^{-\vv_1,p_1}\to\tt W^{-\vv_2,p_2}}\le \bb_t  \exp\big[t\bb_t   \rr_{\vv,p;\dd,k}^t(\mu)^{\theta}\big]   t^{-\frac{i+\xi}2 }.
 \end{equation}
\item[$(2)$] If $\xi:=\vv_1-\vv_2 +\ff{d(p_2-p_1)}{p_1p_2}<1\lor\big(2-i-(\eta-2\kk)^+\big)$, then there exists an increasing function    $\beta: [0,\infty)\to (0,\infty)$  such that    for any  $ t\in (0,\infty) $   and   $  \mu\in \C_{\vv,p; \dd,k}^{t},$
\beq\label{ES11} \|\nn^i  P_{s,t}^\mu\|_{\tt W^{-\vv_1,p_1}\to\tt W^{-\vv_2,p_2}}\le   \bb_t  \exp\big[t\bb_t   \rr_{\vv,p;\dd,k}^t(\mu)^{\theta}\big]  (t-s)^{-\frac{i+\xi}2},
\ \ s\in (0,t). \end{equation}\end{enumerate}
\end{lem}

\begin{proof}
Let $t\in (0,\infty)$ and $   \mu\in \C_{\vv,p; \dd,k}^{t}.$
We will complete the proof  by the following four  steps.

(a)   We first observe that when $p_2>1$,
\beq\label{*}\sup_{0\le r\le s\le t}(s-r)^{\ff i 2}\|\nn^i  {P}_{r,s}^\mu\|_{\tt W^{-\vv_2,p_2} \to \tt W^{-\vv_2,p_2}}  <\infty,\ \ i=0,1.
\end{equation}
By \cite[Proposition 5.4]{HRW25},   \eqref{ZV} implies that for $p_2>1$
\beq\label{ES01} \sup_{0\le r\le s\le t} (s-r)^{\ff i 2}\|\nn^i   P_{r,s}^\mu\|_{\tt L^{p_2} \to \tt L^{p_2}}<\infty, \ \  i=0,1.\end{equation}
Moreover, by the Duhamel formula, see \cite[Proposition 5.5(2)]{HRW25}, we have
\begin{align}\label{DH1}   {P}_{s,t}^\mu f= P_{t-s}^0 f+\int_s^t  {P}_{s,r}^\mu \<b_r(\cdot,\mu_r), \nn   P_{t-r}^0 f\>\d r, \ \  s\in [0,t],  \ f\in \B_b(\R^d).
\end{align}
  We  now prove \eqref{*} by  inducing in $l\in\mathbb N$ for $\vv_2\in [0, lk_0],$  where,  due to \eqref{TJ'},
\beq\label{K0} k_0:= \ff 12 \land [1-(\eta-2\kk)^+]>0.  \end{equation}

Let $\vv_2\in [0, k_0].$ By \eqref{ES01} and   $ \|(1-\DD)^{-\ff{\vv_2}2}\|_{\tt L^{p_2}\to \tt L^{p_2}}<\infty$, we obtain
\beg{align*}& c_1(\mu,t):= \sup_{0\le r\le s\le t} (s-r)^{\ff i 2}  \|\nn^i   P_{r,s}^\mu\|_{\tt L^{p_2}\to \tt W^{-\vv_2, p_2}}\\
&= \sup_{0\le r\le s\le t} (s-r)^{\ff i 2} \|(1-\DD)^{-\ff{\vv_2}2} \nn^i  P_{r,s}^\mu\|_{\tt L^{p_2}\to \tt L^{ p_2}}\\
&\le \|(1-\DD)^{-\ff {\vv_2}2}\|_{\tt L^{p_2}\to \tt L^{p_2}} \sup_{0\le r\le s\le t} (s-r)^{\ff i 2}  \|\nn^i   P_{r,s}^\mu\|_{\tt L^{p_2}\to\tt L^{p_2}}<\infty.\end{align*}
 Combining this with {\bf (A)}, \eqref{Hypin}, \eqref{DH1} and {noting that $s^{\kappa-\frac{\eta}{2}}=s^{-\frac{(\eta-2\kappa)^+}{2}}s^{\frac{(2\kappa-\eta)^+}{2}}$}, we find $c_2(\mu,t)\in (0,\infty)$ increasing in $t\in (0,\infty)$,  such that
  \beg{align*} &\|\nn^i {P}_{r,t}^\mu\|_{\tt W^{-\vv_2,p_2} \to \tt W^{-\vv_2,p_2}}
  \le \|\nn^iP_{t-r}^0\|_{\tt W^{-\vv_2,p_2} \to \tt W^{-\vv_2,p_2}}\\
   &\quad + K_t\int_r^t s^{ \kk } \|\nn^i {P}_{r,s}^\mu\|_{\tt L^{p_2} \to \tt W^{-\vv_2,p_2}}  \|\mu_s\|_{\dd,k*}\|\nn P_{t-s}^0\|_{ \tt W^{-\vv_2,p_2}\to \tt L^{p_2}}\d s\\
 &\le {B_{0}} (t-r)^{-\ff i 2}  + K_tc_1(\mu,t)\rr^t_{\vv,p;\dd,k} (\mu)\int_r^t
 s^{\kk-\ff\eta 2}(s-r)^{-\ff i 2}   (t-s)^{-\frac{\vv_2+1}{2}  }\d s\\
 &\leq {B_{0}}(t-r)^{-\ff i 2}    +  c_2(\mu,t) \int_r^t (s-r)^{-\ff{i+(\eta-2\kk)^+}2} (t-s)^{-\frac{\vv_2+1}{2} } \d s.   \end{align*}
 By \eqref{TJ'}, \eqref{K0}, $\vv_2\in [0,k_0]$ and $i\le 1$,  we have
 $$\frac{\vv_2+1}{2}\lor \ff{i+(\eta-2\kk)^+}  2 <1,\ \ \ \ff{\vv_2+1+(\eta-2\kk)^+  }2\le 1,$$ so that \eqref{LN0} with $\aa=0$ and {$\lambda=0$} implies
 $$ \int_r^t s^{-\frac{i+(\eta-2\kk)^+ }{2}} (t-s)^{-\frac{\vv_2+1}{2} } \d s\le c (t-r)^{-\ff i 2},\ \ r\in [0,t)  $$   for some constant $c\in (0,\infty).$ Therefore,
  \eqref{*}  holds for $\vv_2\in [0,k_0].$

 Assume that for some $l\in \mathbb N$ we have \eqref{*} for $\vv_2\in [0, lk_0]$,  then for $\vv_2\in (lk_0, (l+1)k_0]$,
\beg{align*}& c_l(\mu,t):= \sup_{0\le r\le s\le t} (s-r)^{\ff i 2}\|\nn^i   P_{r,s}^\mu\|_{ \tt W^{-l k_0, p_2}\to \tt W^{-\vv_2, p_2}}\\
 &= \sup_{0\le r\le s\le t} (s-r)^{\ff i 2} \|(1-\DD)^{-\ff{\vv_2-lk_0}2}\nn^i   P_{r,s}^\mu\|_{\tt W^{-l k_0, p_2} \to \tt W^{-lk_0, p_2}} \\
&\le \|(1-\DD)^{-\ff {\vv_2-lk_0}2}\|_{\tt W^{-l k_0, p_2} \to \tt W^{-l k_0, p_2}}\sup_{0\le r\le s\le t} (s-r)^{\ff i 2}\|\nn^i   P_{r,s}^\mu\|_{\tt W^{-l k_0, p_2} \to \tt W^{-l k_0, p_2}}<\infty,\end{align*}
which is increasing in $t\in (0,\infty)$.
 By  combining this with {\bf (A)}, \eqref{Hypin}, \eqref{DH1},
 we find $c_{l+1}(\mu,t)\in (0,\infty)$ increasing in $t\in (0,\infty)$, such that
 \beg{align*} &\|\nn^i {P}_{r,t}^\mu\|_{\tt W^{-\vv_2,p_2}\to\tt W^{-\vv_2,p_2}}
  \le \|\nn^iP_{t-r}^0\|_{\tt W^{-\vv_2,p_2} \to \tt W^{-\vv_2,p_2}} \\
&\qquad + K_t\int_r^t s^{\kk}\|\nn^i {P}_{r,s}^\mu\|_{\tt W^{-l k_0, p_2} \to \tt W^{-\vv_2,p_2}}  \|\mu_s\|_{\dd,k*}\|\nn P_{t-s}^0\|_{ \tt W^{-\vv_2,p_2}\to \tt W^{-lk_0,  p_2}}\d s\\
 &\le {B_{0}} t^{-\ff i 2}  + K_tc_l (\mu,t)\rr^t_{\vv,p;\dd,k} (\mu){B_{\vv_2-lk_0}} \int_r^t s^{\kk-\ff\eta 2} (s-r)^{-\ff i 2}    (t-s)^{-\frac{\vv_2-lk_0+1}{2}  }\d s\\
 &\leq {B_{0}}(t-r)^{-\ff i 2}  +  c_{l+1} (\mu,t) t^{-\ff i 2},\ \ \ r\in [0,t),\  \vv_2\in (l k_0, (l+1)k_0],   \end{align*} where the last step follows  from
 \eqref{LN0} with $\aa=0$ and {$\lambda=0$}, since  $\ff{\eta+i}2-\kk<1$ for $i\le 1$, $\ff{\vv_2-lk_0+1}2<1$
 for $\vv_2-lk_0 \in [0, k_0]$, and  {$\ff{\vv_2-lk_0 +1}2+\frac{(\eta-2\kk)^{+}}{2}\le 1.$}  Hence, \eqref{*}   holds for all $\vv_2\in [0,\infty)$.

 (b) We intend to find  some increasing function $\bb_0: (0,\infty)\to (0,\infty)$ such that
 \beq\label{EE}  \beg{split}&h_{r,t}:= \|\nn^i  P_{r,t}^\mu\|_{\tt W^{-\vv_2,p_2}\to\tt W^{-\vv_2,p_2}}\\
 &\le \bb_0(t) \exp\big[t\bb_0(t)   \rr_{\vv,p;\dd,k}^t(\mu)^{\theta} \big] (t-r)^{-\ff i 2},\
\ \ r\in [0,t), \ i=0,1.\end{split} \end{equation}
  By \eqref{Hypin},   \eqref{DH1} and {\bf (A)}, we obtain that
    \beg{align*} &h_{r,t } \le \|\nn^i P_{r,t}^0\|_{\tt W^{-\vv_2,p_2} \to \tt W^{-\vv_2,p_2}}+ K_t\int_r^t h_{r,s} s^{ \kk}\|\mu_s\|_{\dd,k*} \|\nn P_{t-s}^0\|_{\tt W^{-\vv_2,p_2} \to \tt W^{-\vv_2,p_2}}\d s\\
  &\le   {B_{0}} (t-r)^{-\ff i 2}  +   K_t\rr_{\vv,p;\dd,k}^t (\mu){B_{0}} t^{(\kk-\ff \eta 2)^+}\int_r^t h_{r,s} (s-r)^{-(\ff\eta 2-\kk)^+  } (t-s)^{-\ff 1 2} \d s,\ \ r\in [0,t).\end{align*}
 Combining this with  \eqref{LN0}  for $\aa= \ff 1 \theta=\ff{1-(\eta-2\kk)^+}2$, we find some constant $c(t)\in (0,\infty)$ increasing in $t\in (0,\infty)$, such that
  $$H_{r,t}(\ll):= \sup_{s\in (r,t]} (s-r)^{\ff i 2} h_{r,s} \e^{-\ll (s-r)},\ \ \ s\in (r,t] $$ satisfies
 \beg{align*} H_{r,t}(\ll)&\le {B_{0}} +   K_t\rr_{\vv,p;\dd,k}^t (\mu) {B_{0}}  H_{r,t}(\ll)  t^{(\kk-\ff \eta 2)^+} \\
 &\quad\times \sup_{s\in (r,t]} (s-r)^{\ff i 2}  \int_r^s  (u-r)^{-\ff{(\eta-2\kk)^++i} 2   }    (s-u)^{-\ff 1 2} \e^{-\ll(s-u)} \d u\\
 &\le {B_{0}} +   c(t)     \rr_{\vv,p;\dd,k}^t (\mu) H_{r,t} (\ll)   \ll^{- \ff 1 \theta},\ \ \ r\in [0,t),\ \ll>0.  \end{align*}
  Taking
 $$\ll:= {\big(\frac{1}{2}c(t)   \rr_{\vv,p;\dd,k}^t (\mu)   \big)^{\theta},}$$
 and noting that $H_{r,t}(\ll)<\infty$ for $p_2>1$ due to \eqref{*},  we  find some increasing $\bb_0: (0,\infty)\to (0,\infty)$ such that \eqref{EE} holds for any $p_2 \in (1,\infty]$. Since   $\bb_0(t)$ is uniformly
 in $p_2>1$, by letting $p_2\downarrow 1$ the estimate also holds for $p_2=1.$

 (c) Let $n\in \mathbb N$ such that
 $$\ff {\xi}  n<1,\ \ \ \xi:=   \vv_1-\vv_2 +\ff{d(p_2-p_1)}{p_1p_2},$$
 where we take $n=1$ if $\xi<1$.
Below we prove \eqref{ES1'} by inducing in $0\le l\le n$ such that
\beq\label{67}  \|\nn^i  P_{t}^\mu\|_{\tt W^{-\dd_l,k_l}\to\tt W^{-\vv_2,p_2}}\le  \bb_l(t)  \exp\big[t\bb_l(t)  \rr_{\vv,p;\dd,k}^t(\mu)^{\theta}  \big] t^{-\frac{i}{2}- \ff{l\xi}{2n}}.
 \end{equation}
for some increasing $\bb_l: (0,\infty)\to (0,\infty)$ and
 $$\dd_l:=\vv_2 +\ff l n (\vv_1-\vv_2),\ \ \ k_l:= \ff{n p_1p_2}{(n-l)p_1+lp_2},\ \ \ 0\le l\le n.$$
 In particular, when $l=n$, \eqref{67} reduces to the desired inequality \eqref{ES1'}.

By (b), \eqref{67} holds for $l=0$. Assume that \eqref{67} holds for some $0\le l\le n-1$, where $l=0$ when $\xi<1$, it suffices to verify it for $l+1$ in place of $l$.

To this end, let
\begin{equation*} C_l(\mu,t):=  \bb_l(t) \exp\big[t \bb_l(t)   \rr_{\vv,p;\dd,k}^t(\mu)^{\theta}\big] K_tB_{\vv_1-\vv_2}  \rr_{\vv,p;\dd,k}^t(\mu).\end{equation*}
 By \eqref{Hypin}, \eqref{DH1}, \eqref{67}  and  {\bf (A)}, we obtain
   \beg{align*} &\|\nn^i  {P}_{t}^\mu\|_{\tt W^{-\dd_{l+1},k_{l+1}}\to\tt W^{-\vv_2,p_2}}
  \le \|\nn^i P_{t}^0\|_{ \tt W^{-\dd_{l+1},k_{l+1}}\to\tt W^{-\vv_2,p_2}}\\
  &\qquad + K_t\int_0^t s^{ \kk } \|\nn^i  {P}_{s}^\mu\|_{\tt W^{-\dd_{l},k_{l}}\to\tt W^{-\vv_2,p_2}} \|\mu_s\|_{\dd,k*}\|\nn P_{t-s}^0\|_{ \tt W^{-\dd_{l+1},k_{l+1}}\to\tt W^{-\dd_l,k_l}}\d s\\
 &\le B_{\vv_1-\vv_2}   t^{-\frac{i+(l+1)\xi/n}2 } +   C_l(\mu,t)       \int_0^t  s^{\kk-   \ff{l\xi}{2n} -\ff {i  +\eta}2}  (t-s)^{ -\frac{1+ \xi/n} 2 }\d s.
 \end{align*}
 By $\xi<n$, {\eqref{TJ'}}  and   either $\xi<2-i+2\kk-\eta$ or $\xi<1$ with $l=0$, we have $\ff{1+\xi/n}2<1$ and
 $$\kk-   \ff{l\xi}{2n} -\ff {i  +\eta}2  > -1.$$
 Moreover,  $l+1\le n$ together with  \eqref{TJ'} implies
 $$1+\kk-\ff{l\xi}{2n} -\ff{i+\eta}2 -\ff{1+\xi/n}2 \ge -\ff{i+\xi}2 + 1+\kk -\ff\eta 2 >-\ff{i+\xi}2,$$
so by \eqref{LN0} with $\aa=0$,   we find
      $\bb_{l+1}(t)\in (0,\infty)$  which is increasing in $t>0$ such that \eqref{67} holds for $l+1$ in place of $l$.    Hence, \eqref{ES1'} is proved.

    (d) Finally, let $\xi <2-i-(\eta-2\kk)^+$.  For the above defined $(\dd_l,k_l)_{0\le l\le n}$, we intend to find increasing $\bb_l: (0,\infty)\to (0,\infty) $ such that
    \beq\label{67'}  \|\nn^i  P_{r,t}^\mu\|_{\tt W^{-\dd_l,k_l}\to\tt W^{-\vv_2,p_2}}\le \bb_l(t)  \exp\big[t\bb_l(t)  \rr_{\vv,p;\dd,k}^t(\mu)^{\theta} \big] (t-r)^{-\frac{i}{2}- \ff{l\xi}{2n}  }
    \end{equation} holds for $ 0\le l\le n.$ In particular, when $l=n$ this inequality reduces to the desired \eqref{ES11}.

     By (b), \eqref{67'} holds for $l=0$. Assume that \eqref{67'} holds for some $0\le l\le n-1$, it suffices to verify it for $l+1$ in place of $l$.

 By \eqref{Hypin}, \eqref{DH1}, \eqref{67'}  and  {\bf (A)}, we obtain that for $r\in [0,t)$,
   \beg{align*} &\|\nn^i  {P}_{r,t}^\mu\|_{\tt W^{-\dd_{l+1}, k_{l+1}}\to\tt W^{-\vv_2,p_2}}
  \le \|\nn^i P_{t-r}^0\|_{ \tt W^{-\dd_{l+1},k_{l+1}}\to\tt W^{-\vv_2,p_2}}\\
  &\qquad + K_t\int_r^t s^{ \kk } \|\nn^i  {P}_{r,s}^\mu\|_{\tt W^{-\dd_l,k_l}\to\tt W^{-\vv_2,p_2}} \|\mu_s\|_{\dd,k*}\|\nn P_{t-s}^0\|_{ \tt W^{-\dd_{l+1},k_{l+1}} \to\tt W^{-\dd_l,k_l}}\d s\\
 &\le B_{\vv_1-\vv_2}  (t-r)^{-\frac{i+(l+1)\xi/n}2 } +   {C_l(\mu,t)}  t^{(\kk-\ff\eta 2)^+}    \int_r^t   (s-r)^{-(\ff \eta 2-\kk)^+ -  \ff {i + l\xi/n} 2}  (t-s)^{ -\frac{1+ \xi/n} 2 }\d s.
 \end{align*} 
Since either $\xi< 2-i-(\eta-2\kk)^+ $  with $\xi<n$, or  $\xi<1$ with $l=0$,  we have
 $\ff{\xi/n +1} 2<1$ and
 $$-\Big(\ff \eta 2-\kk\Big)^+ -  \ff {i +l\xi/n}2\ge  -\Big(\ff \eta 2-\kk\Big)^+ -  \ff {i+ \xi} 2>-1.$$
 Moreover,    \eqref{TJ'} together with $l+1\le n$ yields
 $$1- \ff{i+l\xi/n}2 -\Big(\ff\eta 2- \kk\Big)^+-\ff{1+\xi/n}2 >-\ff{i+\xi}2.$$  { So, by \eqref{LN0} with $\aa=0$, }we find
      $\bb_{l+1}(t)\in (0,\infty)$   increasing in $t>0$ such that \eqref{67'} holds for $l+1$ in place of $l$. Then the proof is finished.
\end{proof}

\section{Proof of the existence and uniqueness}

In this section we prove  Theorem \ref{T0}.

Let $X_0$ be $\F_0$-measurable such that $\gg:=\L_{X_0}\in \scr P_{\vv,p*},$  and let $T\in (0,\infty)$ be fixed.
 By Lemma \ref{L2}, for any  $\mu\in \C_{\vv,p;\dd,k}^{\gg,T}$, the SDE
\begin{equation*} \d X_t^\mu=b_t(X_t^\mu,\mu_t)\d t+\d W_t,\ \ X_0^\mu=X_0,\ t\in [0,T]\end{equation*}
has a unique   solution. This provides a map
\begin{equation*} \Phi: \C_{\vv,p;\dd,k}^{\gg,T}\to C^w([0,T];\scr P);\   ( \Phi_t \mu)_{t\in [0,T]} :=(\L_{X_t^\mu})_{t\in [0,T]}.\end{equation*}
So,   the (strong and weak)  well-posedness of $\C_{\vv,p;\dd,k}$-solution for \eqref{E0} up to time $T$, if   $\Phi$ has a unique fixed point in  $\C_{\vv,p;\dd,k}^{\gg,T}.$

For any $n\in\mathbb N$, let $\tau_n(\gg)$ be in \eqref{TT0} for some constant $A_n\in (0,\infty) $ to be determined and
\beq\label{GG1}  \tt\C_{\vv,p; \dd,k}^{\gg, n}:= \Big\{\mu\in \C_{\vv,p;\dd,k}^{\gg,\tau_n(\gg)}:\ \rr_{\vv,p;\dd,k}^{\tau_n(\gg)}(\mu)\le 2B_{\dd-\vv} {\|\gamma\|_{\vv,p*}} \Big\},\end{equation}
where $B_{\dd-\vv}\in (0,\infty)$ is in \eqref{Hypin}. Moreover, for any $\ll\in (0,\infty)$, let
\beq\label{GG2} \hat \C_{\dd,k}^{\gg, n,\ll }:= \Big\{ {\mu\in C^w([0,n],\scr P)}:\ \mu_0=\gg,\ \rr_{\dd,k}^n(\mu):=\sup_{t\in [0,n]} t^{\ff{\dd}2+\ff{d}{2k}} \e^{-\ll t}  \|\mu_t\|_{\dd,k*}\le 2 {B_{\dd}}\Big\}.\end{equation}

 \beg{lem}\label{L3} Assume {\bf(A)} and $ \eqref{TJ'}$ with $p>1$.  Let $B_{\dd-\vv}\in (0,\infty)$ be in $\eqref{Hypin}$, and
 let $\eta:=\dd-\vv+\ff{d(p-k)}{pk}<\kk +\ff 3 2$.  Then the following assertions hold.
 \beg{enumerate} \item[$(1)$]  For any $T\in (0,\infty)$,
 \beq\label{PM}   \Phi: \C_{\vv,p;\dd,k}^{\gg,T}\to \C_{\vv,p;\dd,k}^{\gg,T}.\end{equation}
 \item[$(2)$]  For any $n\in \mathbb N$, there exists a constant $A_n\in (0,\infty)$ such that
 $$\Phi:  \tt\C_{\vv,p; \dd,k}^{\gg, n}\to  \tt\C_{\vv,p; \dd,k}^{\gg, n}.$$  Moreover, if   $\mu $ is a fixed point of
 $\Phi:   \C_{\vv,p;\dd,k}^{\gg,\tau_n(\gg)}\to  \C_{\vv,p;\dd,k}^{\gg,\tau_n(\gg)},$ then   $\mu\in \tt\C_{\vv,p;\dd,k}^{\gg,n}.$
    \item[$(3)$]  If  {$\eta<1+\kappa$} and $ \eqref{TJ'}$ holds for $\vv=0$ and $p=\infty$,   then for any $n\in\mathbb N$ there exists a constant $\ll_n\in (0,\infty)$ such that
    $$\Phi:  \hat\C_{ \dd,k}^{\gg, n,\ll_n}\to   \hat\C_{ \dd,k}^{\gg, n, \ll_n}.$$
    Moreover,  if   $\mu $ is a fixed point of
 $\Phi:   \C_{0,\infty;\dd,k}^{\gg,n}\to   \C_{0,\infty;\dd,k}^{\gg,n},$ then   $\mu\in  \hat\C_{ \dd,k}^{\gg, n,\ll_n}.$
 \end{enumerate}
   \end{lem}

 \beg{proof}      (1) Let $T\in (0,\infty)$ and  $\mu\in \C_{\vv,p;\dd,k}^{\gg,T}. $ We intend to show $\Phi\mu\in \C_{\vv,p;\dd,k}^{\gg,T}.$
 To this end, let
 \beq\label{TA} r=\ff{2}{3+2\kk},\ \ \  \dd':= r\vv+(1-r)\dd,\ \ \ k':= \ff{pk}{(1-r)p+rk}>1.\end{equation}
 By \eqref{TJ'} and  {\ $\eta<\kappa+\frac{3}{2}$ }, we  have
\beq\label{T*} \eta<\ff{3+2\kk}2=  {\ff{2+2\kk}{2-r}=\ff 1 {r}},   \end{equation}
\beq\label{YY}   \dd-\dd' +\ff{d(k'-k)}{k'k}=r\eta,\ \ \ \
   \dd'-\vv +\ff{d(p-k')}{pk'}=(1-r)\eta.  \end{equation}
Then
 {   $ \dd'-\vv +\ff{d(p-k')}{pk'}=(1-r)\eta<2+2\kk-\eta$},
  where $k'>1$,    so that
 \eqref{ES1'} holds for  $i=0$ due to  {by Lemma \ref{L3'}(1), i.e.}
 {          \beq\label{ES1''} \| {P}_t^\mu\|_{\tt W^{-\dd',k'}\to\tt W^{-\vv,p}} \le   \beta_t \exp\big[t\beta_t \rr_{\vv,p;\dd,k}^t(\mu)^{\theta}\big]   t^{-\ff {1-r}2\eta},\ \ t\in (0,T],\ \mu\in \C_{\vv,p;\dd,k}^{\gg,T}.\end{equation}}
 By   \eqref{DH1} for $s=0$, for any $\mu\in \C_{\vv,p;\dd,k}^{\gg,T}$ we have
 \begin{equation*}\beg{split}&(\Phi_t \mu)(f)=\gg( {P}_t^\mu f) \\
 &  = \gg(P_{t}^0  f) + \int_0^t \gg\big( {P}_s^\mu \<b_s(\cdot,\mu_s), \nn   P_{t-s}^0  f\> \big)\d s,\ \ t\in [0, T],\ f\in\B_b(\R^d).\end{split}\end{equation*}
Combining this with  \eqref{Hypin},   \eqref{YY},  \eqref{ES1''}  and {\bf (A)}, we find a constant  $c_1 \in (0,\infty)$ depending on $T$ and $\mu$ such that
\beq\label{Y*}   \beg{split} &\|\Phi_t\mu\|_{\dd,k*}
 \le \|\gamma\|_{\vv,p*}\|P_t^0\|_{ \tt W^{-\dd,k}\to \tt W^{-\vv,p}}\\
 &\qquad + K_t\|\gamma\|_{\vv,p*} \int_0^t s^{  \kk  } \| {P}_s^\mu\|_{\tt W^{-\dd',k'}\to \tt W^{-\vv,p}} \|\mu_s\|_{\dd,k*}\|\nn P_{t-s}^0\|_{  \tt W^{-\dd,k}\to \tt W^{-\dd',k'}}\d s\\
 &\le B_{\dd-\vv}\|\gg\|_{\vv,p*} t^{-\ff\eta 2}  + c_1 \int_0^t s^{ \kk -  \ff{2-r}2 \eta} (t-s)^{-\frac{1}{2}-\ff{r} 2 \eta}\d s,\ \ t\in (0,T].\end{split}\end{equation}
By  \eqref{T*},  {\eqref{TJ'}} and  {\ $\eta<\kappa+\frac{3}{2}$ }, we have
$$  {\kk-\ff{2-r}{2}\eta>-1},\ \ -\ff 1 2 -\ff{r}2\eta>-1,\ \  {\kk- \eta +\ff 1 2\ge -\ff\eta 2}.$$
By \eqref{LN0} for $\aa=0$ and  {$\lambda=0$},  we find  a constant $c\in (0,\infty)$ such that
\beg{align*}  \int_0^t s^{ \kk  - \ff{2-r}2 \eta} (t-s)^{-\frac{1}{2}-\ff{r} 2 \eta}\d s  \le  c t^{-\ff\eta 2},\ \ t>0.\end{align*}
Therefore,  \eqref{Y*}  implies $\Phi\mu\in \C_{\vv,p;\dd,k}^{\gg,T}.$

(2)  Let $\mu\in \tt\C_{\vv,p; \dd,k}^{\gg, n}$ where the constant $A_n$ in  $\tau_n(\gg)$  is to be determined such that
$\Phi\mu\in \tt\C_{\vv,p; \dd,k}^{\gg, n}.$ Then
 {\begin{equation*}\beg{split}& K_nn^{\frac{(2\kappa-\eta)^+}{2}}B_{\delta-\vv}\bb_n\rr^{\tau_n(\gg)}_{\vv,p;\dd,k} (\mu) \exp\Big[n\bb_n\rr_{\vv,p;\dd,k}^{\tau_n(\gg)}(\mu)^{\theta}\Big]\\
&\le D_n(\gg):=  2K_nn^{\frac{(2\kappa-\eta)^+}{2}}B_{\dd-\vv}^2 \|\gg\|_{\vv,p*} \bb_n\exp\Big[n\bb_n (2B_{\dd-\vv}\|\gg\|_{\vv,p*})^{\theta}\Big].\end{split}\end{equation*}}
By combining this with   \eqref{DH1} for $s=0$, {\bf (A)},   \eqref{Hypin}, Lemma \ref{L3'}(1),  we obtain
\beq\label{*D} \beg{split} & \|\Phi_t\mu\|_{\dd,k*} \le \|\gg\|_{\vv,p*} \|P_t^0\|_{\tt W^{-\dd,k}\to\tt W^{-\vv,p}}\\
&\qquad +  {K_t}\int_0^t s^{ \kk  } \|\gg\|_{\vv,p*} \|P_s^\mu\|_{\tt W^{-\dd',k'}\to \tt W^{-\vv,p}}   \|\mu_s\|_{\dd,k*} \|\nn P_{t-s}^0 \|_{\tt W^{-\dd, k}\to\tt W^{-\dd',k'}}\d s\\
& {\le B_{\dd-\vv} \|\gg\|_{\vv,p*}t^{-\ff\eta 2} + \|\gg\|_{\vv,p*}  2K_nB_{\dd-\vv}^2 \|\gg\|_{\vv,p*} \bb_n\exp\Big[n\bb_n (2B_{\dd-\vv}\|\gg\|_{\vv,p*})^{\theta}\Big]}\\
&\qquad\qquad\qquad\qquad\quad\quad\quad {\times \int_0^t
s^{\kk-\ff{2-r}2\eta} (t-s)^{-\ff 1 2 -\ff{r}2 \eta}\d s}\\
&\le  B_{\dd-\vv} \|\gg\|_{\vv,p*}t^{-\ff\eta 2}+ \|\gg\|_{\vv,p*}  D_n(\gg) c t^{\ff 1 \theta-\frac{\eta}{2}},   \ \ t\in (0,\tau_n(\gg)],\end{split}\end{equation}
where $c\in (0,\infty)$ is a constant due to \eqref{LN0} for $\aa=0$ and  {$\lambda=0$}.
Taking $A_n\in (0,\infty)$   such that
\beq\label{RB} \Big(\ff{B_{\dd-\vv} }{D_n(\gg)c}\Big)^{\theta}\ge \Big(A_n \e^{A_n \|\gg\|_{\vv,p*}^{\theta}}\Big)^{-1}, \end{equation}
by the definition of $\tau_n(\gg) $ we derive
$$\|\Phi_t\mu\|_{\dd,k*}\le 2 B_{\dd-\vv}  {t^{-\frac{\eta}{2}}}\|\gg\|_{\vv,p*},\ \ \ t\le \tau_n(\gg).$$
  { So,  $\Phi\mu\in\tt\C_{\vv,p; \dd,k}^{\gg, n}.$ }

If $\mu\in \C_{\vv,p; \dd,k}^{\gg, T}$ is a fixed point of $ \Phi$, then $ \Phi\mu=\mu$.
  { Noting that $\rr_{\vv,p;\dd,k}^{t}(\mu)$ is non-decreasing in $t$,
\begin{equation*} \rr_{\vv,p;\dd,k}^{t+}(\mu):= \lim_{\vv\downarrow 0} \rr_{\vv,p;\dd,k}^{(t+\vv)\land T}(\mu),\ \
\rr_{\vv,p;\dd,k}^{t-}(\mu):= \lim_{\vv\downarrow 0}  \rr_{\vv,p;\dd,k}^{(t-\vv)^+}(\mu)\end{equation*}
exist and are non-decreasing for $t\in (0,T].$}
So, the first inequality in \eqref{*D}, \eqref{Hypin} and Lemma \ref{L3'}(1) yield
\beg{align*} &\|\mu_t\|_{\dd,k*}= \|\Phi_t\mu\|_{\dd,k*}\le B_{\dd-\vv}\|\gg\|_{\vv,p*}t^{-\ff\eta 2}\\
&\quad + { K_tB_{\delta-\vv}\|\gg\|_{\vv,p*} \bb_t \exp\big[t\bb_t \rr_{\vv,p;\dd,k}^{t-} (\mu)^{\theta}\big] \rr_{\vv,p;\dd,k}^{t-}(\mu) \int_0^t s^{\kk-\ff{2-r}{2}\eta}(t-s)^{-\ff 1 2 -\ff{r}2 \eta}\d s}\\
&\le B_{\dd-\vv}\|\gg\|_{\vv,p*}t^{-\ff\eta 2}+   { K_tB_{\delta-\vv}\|\gg\|_{\vv,p*} t^{\frac{(2\kappa-\eta)^+}{2}}\bb_t \exp\big[t\bb_t \rr_{\vv,p;\dd,k}^{t-} (\mu)^{\theta}\big] \rr_{\vv,p;\dd,k}^{t-} (\mu) ct^{\ff 1 \theta-\ff\eta 2}}.\end{align*}
Thus,
 \beq\label{PI} \beg{split}&\rr^{t+}_{\vv,p;\dd,k}(\mu)\le B_{\dd-\vv}  \|\gg\|_{\vv,p*} \\
 &\quad +  { K_tB_{\delta-\vv}\|\gg\|_{\vv,p*} t^{\frac{(2\kappa-\eta)^+}{2}} \bb_n \exp\big[n\bb_n \rr_{\vv,p;\dd,k}^{t-} (\mu)^{\theta}\big] \rr_{\vv,p;\dd,k}^{t-} (\mu) ct^{\ff 1 \theta}},\ \ t\in (0, \tau_n(\gg)]. \end{split}\end{equation}
Then
\begin{equation*} \rr^{0+}_{\vv,p;\dd,k}(\mu)  \le  B_{\dd-\vv} \|\gg\|_{\vv,p*}.\end{equation*}  This and the right continuity of $\rr_{\vv,p;\dd,k}^{t+}(\mu)$ in $t\ge 0$ imply
$$s_0:=  \tau_n(\gg)\land \inf\big\{t\in (0, \tau_n(\gg)]:\  \rr^{t+}_{\vv,p;\dd,k}(\mu)\ge 2 B_{\dd-\vv}  \|\gg\|_{\vv,p*}\big\}>0,$$
where $\inf\emptyset:=\infty$ by convention.
 If $s_0< \tau_n(\gg)$,    by the non-decreasing of $\rr^{t}_{\vv,p;\dd,k}$ in $t\ge 0$,   we obtain
 $$  \rr^{s_0+}_{\vv,p;\dd,k}(\mu)\ge 2 B_{\dd-\vv}  \|\gg\|_{\vv,p*}\ge \rr^{s_0-}_{\vv,p;\dd,k}(\mu),$$
 so that  \eqref{PI} yields
 \beg{align*}& 2 B_{\dd-\vv}    \|\gg\|_{\vv,p*} \le  \rr^{s_0+}_{\vv,p;\dd,k}(\mu)  \leq B_{\dd-\vv}  \|\gg\|_{\vv,p*}  +  D_n(\gg)c  \|\gg\|_{\vv,p*}    s_0^{\theta}.\end{align*}
So, \eqref{RB} and the definition of $\tau_n(\gg)$ imply
$$ s_0\ge \Big(\ff{B_{\dd-\vv}}{D_n(\gg)c}\Big)^{  \theta} \ge   \tau_n(\gg),$$
which  contradicts  to $s_0< \tau_n(\gg)$. Hence, $s_0\ge \tau_n(\gg)$, which implies
    $\mu\in\tt\C_{\vv,p;\dd,k}^{\gg, n}.$

    (3) Let $\mu\in \hat\C_{\dd,k}^{\gg, n,\ll}$, where $\gg\in \scr P, n\in \mathbb N$ and $\ll\in (0,\infty)$.
  {   By \eqref{DH1}, \eqref{Hypin}, {\bf(A)} and noting that $\|\gg\|_{0,\infty*}= 1$ for $\gg\in \scr P$,} we obtain
 {   \beg{align*}  \|\Phi_t\mu\|_{\dd,k*} &\le B_{\dd}t^{-\ff{\eta}2} +K_n B_{0}\int_0^t s^\kk \|\Phi_s\mu\|_{\dd, k*} \|\mu_s\|_{\dd,k*} (t-s)^{-\ff 1 2}\d s\\
    &\le B_{\dd} t^{-\ff{\eta}2} + 2B_{0}B_\delta K_n\int_0^t s^{\kk-\ff\eta 2}  \|\Phi_s\mu\|_{\dd,k*} (t-s)^{-\ff 1 2}\d s,\ \ t\in (0,n].\end{align*}}
Combining this with  \eqref{LN0} for $\alpha=\ff{1-(\eta-2\kk)^+}2$, when  $\eta<1+\kappa$ such that $\kk-\eta >-1,$ we find  a constant $D_n\in (0,\infty)$ such that
    $$H_{n,\ll}:= \sup_{t\in (0,n]} t^{\ff\eta 2} \|\Phi_t\mu\|_{\dd, k*} \e^{-\ll t}$$
    satisfies
    \beg{align*} H_{n,\ll}&\le B_{\dd}+  {2 K_nB_0B_{\dd}}  H_{n,\ll} \sup_{t\in (0,n]} t^{\ff\eta 2} \int_0^t s^{\kk- \eta  }   (t-s)^{-\ff 1 2} \e^{-\ll (t-s)} \d s\\
    &\le B_{\dd} +  D_n  H_{n,\ll} \ll^{-\ff{1-(\eta-2\kk)^+}2},\ \ \ \ll>0.\end{align*}
    Noting that  \eqref{PM} for $\vv=0$ and $p=\infty$ implies $H_{n,\ll}<\infty$, taking
    $$\ll_n:= (2D_n)^{-\ff 2 {1-(\eta-2\kk)^+}},$$
    we derive $H_{n,\ll_n}\le 2 B_{\dd}$, so that $\Phi: \hat\C_{ \dd,k}^{\gg, n,\ll_n}\to   \hat\C_{\dd,k}^{\gg, n, \ll_n}.  $

    If $\mu\in \C_{0,\infty;\dd,k}^{\gg,n}$ is a fixed point of $\Phi$, then $\mu\in \hat\C_{\dd,k}^{\gg, n,\ll_n}$ can be proved by the same argument as in step (2) to verify that
    $$s_0:= n\land \inf\big\{t\in (0, n]:\  \rr^{t+}_{0,\infty;\dd,k}(\mu)\ge  {2 B_{\dd}}  \big\}=n.$$

 \end{proof}

Let $r$ be in \eqref{TA},  by   \eqref{T*}  {and $\eta r<1$},  we have
\beg{align*}
  \frac{1}{\theta_1}:=\ff 1 2 -\ff{r}2\eta -\Big(\ff {1-r} 2 \eta-\kk\Big)^+ =\Big(\ff 1 2 - \ff{r} 2\eta\Big)\land \ff 1\theta >0.\end{align*}
We have the following estimate on the Lipschitz  continuity of $\Phi$ under $\rr_{\vv,p;\dd,k}^{\ll,T}.$

\beg{lem}\label{NL} Assume {\bf (A)} and $\eqref{TJ'}$ with $p>1$. Then there exists    increasing $C: (0,\infty)\to (0,\infty)$ such that
\beg{align*} &\rr_{\vv,p;\dd,k}^{\ll,T}(\Phi\mu,\Phi\nu)\le \ll^{-\ff 1 {\theta_1}} \|\gg\|_{\vv,p*}  C(T)  \exp\Big[ C(T)\big(\rr_{\vv,p;\dd,k}^T(\mu)+\rr_{\vv,p;\dd,k}^T(\nu)\big)^{ \theta}\Big] \rr_{\vv,p;\dd,k}^{\ll,T}(\mu,\nu),\\
&\qquad  \ T\in (0,\infty),\ \mu,\nu\in \C_{\vv,p;\dd,k}^{\gg,T},\ \ll\ge \big(C(T) \rr_{\vv,p;\dd,k}^T(\mu)\big)^\theta.\end{align*}
 \end{lem}

\beg{proof}  Let $r,\dd'$ and $k'$ be in \eqref{TA}.  By \eqref{DH1} for $s=0$, {\bf (A)}, \eqref{Hypin} and \eqref{ES1''}, we obtain
 {\beg{align*} &\|P_t^\mu-P_t^\nu\|_{\tt W^{-\dd',k'}\to\tt W^{-\vv,p}}=\sup_{f\in\B_b(\R^d), \|f\|_{\tt W^{-\dd',k'}}\le 1} \|P_t^\mu f-P_t^\nu f\|_{\tt  W^{-\vv,p}}\\
&\le \int_0^t \sup_{f\in\B_b(\R^d), \|f\|_{\tt W^{-\dd',k'}}\le 1}\Big[\big\|(P_s^\mu -P_s^\nu) \<b_s(\cdot,\mu_s),\nn P_{t-s}^0f\>\big\|_{\tt  W^{-\vv,p}}\\
&\qquad \qquad\qquad\qquad \qquad\qquad +
\big\|P_s^\nu \<b_s(\cdot,\mu_s)-b_s(\cdot,\nu_s),\nn P_{t-s}^0 f\>\big\|_{\tt  W^{-\vv,p}}\Big]\d s\\
&\le B_{0}K_t \rr_{\vv,p;\dd,k}^t(\mu) \int_0^t s^{\kk-\ff\eta 2} \|P_s^\mu-P_s^\nu\|_{\tt W^{-\dd',k'}\to\tt W^{-\vv,p}} (t-s)^{-\ff 1 2 }\d s \\
&\quad +B_{0}K_tc_0(t) \exp\big[c_0(t)  \rr_{\vv,p;\dd,k}^{t}(\nu)^{\theta}\big]  \rr_{\vv,p;\dd,k}^{\ll,t} (\mu,\nu) \int_0^t s^{\kk-\ff{2-r}2 \eta }\e^{\ll s} (t-s)^{-\ff 1 2}\d s.\end{align*}}
Combining this with \eqref{LN0} for $\aa= \ff 1 \theta$ and $\aa=0$ respectively,  {$\aa_1= (\ff{2-r}2\eta-\kk)^+<1$} due to \eqref{T*}, and $\aa_2=\ff 1 2$, we find   $c_1(t)\in (0,\infty)$ increasing in $t>0$ such that
 {\beg{align*} &H_t(\ll):=\sup_{u\in (0,t]}\e^{-\lambda u}u^{\ff {1-r} 2 \eta} \|P_u^\mu-P_u^\nu\|_{\tt W^{-\dd',k'}\to \tt W^{-\vv,p}}\\
&\le B_{0}K_t \Big[ \rr_{\vv,p;\dd,k}^t(\mu) H_t(\ll)+c_0(t) \exp\big[c_0(t)  \rr_{\vv,p;\dd,k}^{t}(\nu)^{\theta}\big]  \rr_{\vv,p;\dd,k}^{\ll,t} (\mu,\nu) \Big]\\
&\qquad \times \sup_{u\in (0,t]} u^{\ff {1-r} 2 \eta+(\kk-\ff{2-r}2\eta)^+} \int_0^u s^{-(\ff{2-r}2 \eta-\kk)^+} (u-s)^{-\ff 1 2} \e^{-\ll(u-s)}\d s \\
&\le c_1(t)  \rr_{\vv,p;\dd,k}^t(\mu) H_t(\ll)\ll^{-\ff 1 \theta} + c_1(t)  \exp\big[c_0(t)  \rr_{\vv,p;\dd,k}^{t}(\nu)^{\theta}\big]  \rr_{\vv,p;\dd,k}^{\ll,t} (\mu,\nu)\end{align*}}
holds for all $\ll>0$. By \eqref{ES1'} for $i=0$,  {$(\vv_1,p_1)=(\dd',k')$ and $(\vv_2,p_2)=(\vv,p)$,} we have $H_T(\ll)<\infty$ for $\ll>0$. So, taking
\beq\label{BLL} \bar{\lambda}(T,\mu)= \big(2c_1(T)  \rr_{\vv,p;\dd,k}^T(\mu)\big)^{\theta},\end{equation}
  we find  a constant $c_2(t)\in (0,\infty)$ increasing in $t$ such that  {for any $\lambda\ge \bar{\lambda}(T,\mu)$,}
 {\beq\label{TL}\beg{split} H_t(\lambda)
&\le  c_2(t)  \exp\Big[c_2(t)  \big(\rr_{\vv,p;\dd,k}^{t}(\nu)+\rr_{\vv,p;\dd,k}^t(\mu)\big)^{\theta} \Big]  \rr_{\vv,p;\dd,k}^{\ll,t} (\mu,\nu),\\
&\qquad\qquad \ \ \ t\in (0,T],\ \mu,\nu\in \scr C_{\vv,p;\dd,k}^{\gg,t}.\end{split} \end{equation}}
Let $(r,\dd',k') $ be in \eqref{TA}. By {\bf(A)} and \eqref{DH1} for $s=0$, we have
\beg{align*} &\|\Phi_t\mu-\Phi_t\nu\|_{\dd,k*}= \sup_{f\in \B_b(\R^d), \|f\|_{\tt W^{-\dd,k} }\le 1} \big|\gg(P_t^\mu-P_t^\nu)f\big|\\
&\le \|\gg\|_{\vv,p*} \int_0^t  \sup_{f\in \B_b(\R^d), \|f\|_{\tt W^{-\dd,k}} \le 1}  \big\|P_t^\mu\<b_s(\cdot,\mu_s),\nn P_{t-s}^0f\> - P_t^\nu\<b_s(\cdot,\nu_s),\nn P_{t-s}^0f\>\big\|_{\tt W^{-\vv,p}} \d s\\
&\le \|\gg\|_{\vv,p*} K_t \int_0^t s^\kk \|P_s^\mu-P_s^\nu\|_{\tt W^{-\dd',k'}\to\tt W^{-\vv,p}}\|\mu_s\|_{\dd,k*} \|\nn P_{t-s}^0\|_{\tt W^{-\dd,k}\to \tt W^{-\dd',k'}}\d s\\
&\quad + \|\gg\|_{\vv,p*} K_t \int_0^t s^\kk \| P_s^\nu\|_{\tt W^{-\dd',k'}\to\tt W^{-\vv,p}}\|\mu_s-\nu_s\|_{\dd,k*} \|\nn P_{t-s}^0\|_{\tt W^{-\dd,k}\to \tt W^{-\dd',k'}}\d s.\end{align*}
Combining this with \eqref{Hypin}, \eqref{YY} and \eqref{ES1''}  we find $c_3(t)\in (0,\infty)$ increasing in $t$ such that
\begin{equation*}\beg{split} & \|\Phi_t\mu-\Phi_t\nu\|_{\dd,k*}\\
&\le   c_3(t)\|\gg\|_{\vv,p*} \rr_{\vv,p;\dd,k}^t(\mu) \int_0^t  s^{\kk-\ff\eta 2} \|P_s^\mu-P_s^\nu\|_{\tt W^{-\dd',k'}\to\tt W^{-\vv,p}} (t-s)^{-\ff1 2-\ff{r}2\eta}\d s\\
& + c_3(t)\|\gg\|_{\vv,p*} \exp\big[c_0(t)\rr_{\vv,p;\dd,k}^t(\mu)^\theta\big] \int_0^t  \|\mu_s-\nu_s\|_{\dd,k*}s^{\kk -\ff {1-r}2 \eta}  (t-s)^{-\ff1 2-\ff{r}2\eta}\d s.\end{split} \end{equation*}
Thus, by 
 {\eqref{TL},} and noting that $1+2\kk-\eta\ge 0$ due to \eqref{TJ'}, we find a constant $c_4(t)\in (0,\infty)$ increasing in $t$ such that
\beq\label{YM} \beg{split} &\rr_{\vv,p;\dd,k}^{\ll,T}(\Phi\mu,\Phi\nu)= \sup_{t\in (0,T]} t^{\ff\eta 2}\e^{-\ll t}   \|\Phi_t\mu-\Phi_t\nu\|_{\dd,k*}\\
&\le  \|\gg\|_{\vv,p*} c_3(T) \exp\Big[c_4(T) \big(\rr_{\vv,p;\dd,k}^{T}(\nu)+\rr_{\vv,p;\dd,k}^T(\mu)\big)^{\theta} \Big]  \rr_{\vv,p;\dd,k}^{\ll,T} (\mu,\nu)\\
&\qquad \times\sup_{t\in (0,T]}  {t^{\frac{\eta}{2}}} t^{(\kk -\ff{2-r}2\eta)^+}
\int_0^t   s^{  -(\ff{2-r}2\eta-\kk)^+}   (t-s)^{-\ff 1 2-\ff{r}2\eta}\e^{-\ll (t-s)}\d s,\\
& {\qquad \qquad \ \ \ t\in (0,T],\   \mu,\nu\in \C_{\vv,p;\dd,k}^{\gg,T},\ \ll\ge \bar\ll(T,\mu).} \end{split}\end{equation}
By \eqref{LN0} for $\aa=\theta_1$, $\aa_1=(\ff{2-r}2\eta-\kk)^+<1$  and $\aa_2= \ff 1 2+\ff{r}2\eta<1$ due to \eqref{T*},  we find a constant $c_5(T)\in (0,\infty)$ increasing in $T$ such that
\beg{align*} & \sup_{t\in (0,T]} t^{\ff\eta 2+ (\kk -\ff{2-r}2\eta)^+}  \int_0^t   s^{  -(\ff{2-r}2\eta-\kk)^+}  (t-s)^{-\ff 1 2-\ff{r}2\eta}\e^{-\ll(t-s)}\d s \\
& {\le c_5(T) \ll^{-\theta_1},\ \ \ \ll\ge \bar{\lambda}(T,\mu),\  t\in (0,T].}\end{align*}
Combining this with \eqref{BLL} and \eqref{YM}, we derive   the desired estimate  for some $C(T)\in (0,\infty)$ increasing in $T$.   \end{proof}

We are now ready to prove the following result, which implies Theorem \ref{T0}.

 \beg{prp}\label{PP} Assume {\bf (A)} and  let $\vv\in [0,\dd]$ and $p\in [k,\infty] $ satisfying $\eqref{TJ'}$.
 \beg{enumerate} \item[$(1)$]
 For any $\F_0$-measurable  initial value $X_0$ with $\gg=\L_{X_0}\in \scr P_{\vv,p*}$, $\eqref{E0}$ has a unique maximal  (weak and strong)
  $\C_{\vv,p;\dd,k}$-solution.
 \item[$(2)$] For any $n\in \mathbb N$, there exist $A_n,\ll_n\in (0,\infty)$ such that $\eqref{TT0}$ and $\eqref{EST}$ hold.
 \item[$(3)$] There exists an increasing function $C_\gg: [1,\infty)\to (0,\infty)$, which does not depend on $\gg$ when $\vv=0,p=\infty$ and $\eta<1$, such that $\eqref{NES}$ holds.
 \end{enumerate}\end{prp}

\beg{proof}   According to the proof of \cite[Theorem 2.1]{HRW25}, the first assrtion follows from the second. So, it suffices to show that for any $n\in \mathbb N$ and initial value $X_0$ with $\gg\in \L_{X_0}\in \scr P_{\vv,p*},$
$\eqref{E0}$ has a unique weak/strong $\C_{\vv,p;\dd,k}$-solution up to $\tau_n(\gg)$ such that
\beq\label{ERR} \E\bigg[\sup_{t\in (0,\tau_n(\gg)]} |X_t|^q\bigg|\F_0\bigg] \le c_{n,q}(\gg) \big(1+|X_0|^q\big),\ \ q\in (0,\infty),\end{equation}
where $c_{n,q}(\gg)\in (0,\infty)$, which does not depend on $\gg$ when    $\vv=0$ and $p=\infty$.
Below, we simply denote $ \tau_n=\tau_n(\gg).$

(a)  By Lemma \ref{L2} and \eqref{PM}, $\eqref{E0}$ has a unique weak/strong $\C_{\vv,p;\dd,k}$-solution up to $\tau_n$ provided
$\Phi:  \C_{\vv,p;\dd,k}^{\gg,\tau_n}\to \C_{\vv,p;\dd,k}^{\gg,\tau_n}$ has a unique fixed point.

By Lemma \ref{L3},
  all fixed points in $\C_{\vv,p;\dd,k}^{\gg, \tau_n}$ of $ \Phi$ are included in $\tt\C_{\vv,p;\dd,k}^{\gg, n}$ when $\vv>0$ or $p<\infty$, and in $\hat \C_{\dd,k}^{\gg,n,\ll_n}$ when $\vv=0$ and $p=\infty$.
So, \eqref{EST} holds for any (weak) $\C_{\vv,p;\dd,k}$-solutions of $\eqref{E0}$ with initial distribution $\gg$ up to time $ \tau_n$.
By Lemma \ref{L1} and the contractive  fixed point theorem, it suffices to find $\ll\in (0,\infty)$ such  that the map
$  \Phi$ is contractive under the metric  $\rr^{\ll,\tau_n}_{\vv,p;\dd,k}$ on $ \tt\C_{\vv,p; \dd,k}^{\gg, n}$ when $\vv>0$ or $p<\infty$, and on $ \hat \C_{\dd,k}^{\gg, n,\ll_n}$ otherwise.

 Let  $\vv>0$ or $p<\infty$.
  By Lemma \ref{NL}  and \eqref{GG1},  we can find $C_{n,\gg}\in (0,\infty)$ such that
$$\rr_{\vv,p;\dd,k}^{\ll,\tau_n}(\Phi\mu,\Phi\nu)\le \ll^{-\ff 1 {\theta_1}} C_{n,\gg} \rr_{\vv,p;\dd,k}^{\ll,\tau_n}(\mu,\nu), \ \ \mu,\nu\in  \tt\C_{\vv,p; \dd,k}^{\gg, n},\ \ \lambda\ge  C_{n,\gg}.$$
Similarly, when $\vv=0$ and $p=\infty$, by Lemma \ref{NL}  and \eqref{GG2},  we find a constant $C_n\in (0,\infty)$ uniformly in $\gg\in \scr P$ such that
$$\rr_{\vv,p;\dd,k}^{\ll,n}(\Phi\mu,\Phi\nu)\le \ll^{-\ff 1 {\theta_1}} C_{n} \rr_{\vv,p;\dd,k}^{\ll,n}(\mu,\nu),\ \  \mu,\nu\in  \hat\C_{\dd,k}^{\gg, n,\ll_n},\lambda\ge C_n.$$
So, in any case $\Phi$ is contractive under the metric  $\rr^{\ll, \tau_n}_{\vv,p;\dd,k}$ when
$\ll>0$ is large enough.  So, $\eqref{E0}$ has a unique weak/strong $\C_{\vv,p;\dd,k}$-solution up to $\tau_n$.

(b) Let $X_t$ be the unique solution up to time $\tau_n$. Then
$$(\mu_t=\L_{X_t})_{t\in [0,\tau_n]}\in  \beg{cases} \tt\C_{\vv,p; \dd,k}^{\gg, n}, &\text{if}\ \vv>0\ \text{or}\ p<\infty,\\
 \hat\C_{\dd,k}^{\gg, n,\ll_n}, &\text{if}\ \vv=0, p=\infty.\end{cases}$$
Then by \eqref{TJ'} and \eqref{EST},  there exists $q'>2$   such that $\|b^\mu\|_{\tt L_{q'}^\infty(\tau_n)}$ in   \eqref{ZV}  is bounded above by a constant
$D_{n,\gg}\in (0,\infty)$ depending on $n$ and $\|\gg\|_{\vv,p*}$. So, \eqref{ERR} follows from \cite[Theorem 1.3.1]{RW24}, where $c_{n,q}(\gg)$ is uniform in $\gg\in \scr P$ when $\vv=0$ and $p=\infty$,
since   $\|\gg\|_{0,\infty*}=1$.
   \end{proof}

\section{Proof of regularity estimates }


 By Theorem \ref{T0},
for any  $\gamma\in\scr P_{\vv,p*}$ and $T\in (0,\tau(\gg))$,   we have
$ (P_t^*\gg)_{t\in [0,T]}\in \C_{\vv,p;\dd,k}^{\gg,T}.$
Let $P_{s,t}^\gg= {P}_{s,t}^\mu$ be defined in \eqref{SM} for   $\mu_t=P_t^\ast\gamma$, i.e.
\begin{align*}  P_{s,t}^\gamma f(x)=\E [f({X}_{s,t}^{\gg,x})],\ \ 0\le s\le t< \tau(\gamma),\ f\in\scr B_b(\R^d),\ x\in\R^d,
 \end{align*}
where for fixed   $(s,x)\in [0,\tau(\gg))\times\R^d$,    $({X}_{s,t}^{\gg,x})_{  t\in [s, \tau(\gamma))}$ is the unique solution to  the SDE
$$\d  {X}_{s,t}^{\gg,x}=b_t({X}_{s,t}^{\gg,x}, P_t^\ast\gamma)\d t+\d W_t,\ \ {X}_{s,s}^{\gg,x}=x,\ t\in [s,\tau(\gg)).$$
Simply denote $P_t^\gg=P_{0,t}^\gg$ for $t\in [0,\tau(\gamma))$.

 We now  prove that \eqref{TJ} implies
 \beq\label{JJ} 1+(2\kk-\eta)^+-\eta> \Big(\vv+\ff d p-\ff d k\Big)^+,\end{equation}
 which is claimed in Theorem \ref{T02}.

 If $\eta<2\kk$, then $\eta<1\lor (\ff 1 2+\kk)$  implies
 $$1+(2\kk-\eta)^+-\eta= 1+2\kk-2\eta>0,$$
and by $\dd<1+(2\kk-\eta)^+= 1+2\kk-\eta$,
$$1+(2\kk-\eta)^+-\eta- \Big(\vv+\ff d p-\ff d k\Big)= 1+2\kk-2\eta+\eta-\dd= 1+2\kk-\eta-\dd >0,$$ so that
\eqref{JJ} holds.

If $\eta\ge 2\kk$ then \eqref{TJ} implies $\kk<\ff 1 2$ and $\dd<1$, so that
$$1+(2\kk-\eta)^+-\eta= 1-\eta  >0,$$
and
 $$1+(2\kk-\eta)^+-\eta- \Big(\vv+\ff d p-\ff d k\Big) = 1-\eta+\eta-\dd =1-\dd>0.$$
 Hence, \eqref{TJ} implies \eqref{JJ}.

\beg{proof}[\textbf{Proof of Theorem \ref{T02}(1)}] All constants   $\{c_i(t)\}_{i\ge 1}\subset (0,\infty)$ below are  increasing in  $t>0$. For any $\gg\in \scr P_{\vv,p*},$ define $P_t^{\gg*}:\scr P\to\scr P$ by
\begin{equation*} (P_t^{\gamma\ast}\nu)(f):=\int_{\R^d}\big(P_t^\gamma f(x)\big)\nu(\d x),\ \ f\in\scr B_b(\R^d), \ t\in [0,\tau(\gg)),\ \nu\in \scr P.\end{equation*}
Then
$ P_t^*\gg= P_t^{\gg*}\gg$  for $\gg\in \scr P_{\vv,p*},$  so that
\beq\label{EQ}\|P_t^*\gg-P_t^*\tt\gg\|_{\dd,k*}\le  \|P_t^{\tt\gg*}\gg-P_t^{\tt\gg*}\tt\gg\|_{\dd,k*}
+ {\|P_t^{\gg*}\gg-P_t^{\tt\gg*}\gg\|_{\dd,k*}}.\end{equation}
In the following, we estimate these two terms respectively for fixed $\gg,\tt\gg\in \scr P_{\vv,p*}$.

 (a) Let $\pi\in \C(\gg,\tt\gg)$ such that
\beq\label{CPL} \W_q(\gg,\tt\gg)= \bigg(\int_{\R^d\times\R^d}|x-y|^q\pi(\d x,\d y)\bigg)^{\ff 1 q}.\end{equation}
By $\vv\in [0,\ff{d(p-1)}p]$, we have $p_0:= \ff{dp}{d+\vv p}\in [1,p]$ and
$$\ff 1 {p_0}= \ff 1 p+ \ff\vv d.$$
Then   the Sobolev embedding theorem implies
$$\|\cdot\|_{\tt W^{-\vv,p}}\le c_0 \|\cdot\|_{\tt L^{p_0}}$$
  for some constant $c_0\in (0,\infty),$ so that
\beq\label{SB} \|\gg\|_{p_0*}:=\|\gg\|_{0,p_0*} \le  {c_0 \|\gg\|_{\vv,p*}},\ \ \gg\in \scr P_{\vv,p*}.\end{equation}
By the upper bound condition in \eqref{QY}, we have
$$k_0:=\ff{p_0q}{q-1}\ge k.$$
 By  the definitions of $\eta, p_0$ and $k_0$, we obtain
$$\xi(q):=\dd+\ff{dpq-{k(d+\vv p)(q-1)}}{pqk}=\dd+\ff{d(k_0-k)}{k_0k}
=\eta +\ff 1 q \Big(\vv+ \ff d p\Big).$$
Then \eqref{QY} implies \eqref{XQ}, so that
$$\eta\le \xi(q) <1+(2\kk-\eta)^+\le 1\lor \big(2-i+2\kk-\eta\big),\ \ i=0,1.$$
Moreover, by the second inequality in \eqref{TJ}, we conclude $\dd+\frac{d}{k}< 1\vee(2+2\kk-\eta)$.
Thus, by Lemma \ref{L3'}(1), there exists increasing $\bb: (0,\infty)\to (0,\infty)$ such that
\beq\label{X1} \begin{split}&\|\nn^i P_t^{\tt\gg}\|_{\tt W^{-\dd,k}\to \tt W^{0,k_0}} \le \bb_t \exp\big[t\bb_t k_t(\tt\gg)^\theta\big] t^{-\ff{i+\xi(q)}2},\ \ \ i=0,1,\ t\in (0,\tau(\tilde{\gg})),\\
&\|P_t^{\tt\gg} \|_{\tilde{W}^{-\delta,k}\to \tilde{W}^{0,\infty}}\leq \bb_t \exp\big[t\bb_t k_t(\tt\gg)^\theta\big] t^{-\ff{\delta+\frac{d}{k}}2}, \ \ t\in (0,\tau(\tilde{\gg})).
\end{split}\end{equation}
Next,  consider the maximal functional
\begin{equation*}\scr M f(x):= \sup_{r\in (0,1)}\ff 1 {|B(x,r)|} \int_{B(x,r)} f(y)\d y,\ \ x\in\R^d,\end{equation*}
   for a nonnegative measurable function $f$.
By \cite[Lemma 2.1]{XXZZ}  and  $P_t^\gg f\in C(\R^d)$  for $f\in \B_b(\R^d) $ due to \eqref{X1},  we find a constant $k_1\in (0,\infty)$ such that
  \beg{align*} &|P_t^{\tt\gg} f(x)-P_t^{\tt\gg} f(y)|\le k_1 |x-y|\big(\scr M|\nn P_t^{\tt\gg} f|(x)+\scr M|\nn P_t^{\tt\gg} f|(y)+ {\|P_t^{\tt\gg} f\|_\infty}\big),\\
 & \big\|\scr M|\nn P_t^{\tt\gg} f|\big\|_{\tt L^{k_0}}\le k_1 \|\nn P_t^{\tt\gg} f\|_{\tt L^{k_0}}, \ \  \ t\in (0,\tau(\gg)),\ x,y\in\R^d.
\end{align*}
  Combining this with H\"older's inequality, \eqref{CPL}, \eqref{SB}, \eqref{X1}  {and $\|\gamma\|_{\vv,p*}\geq 1$}, we find    $k_2, c_1(t)\in (0,\infty)$ such that for any $f\in \B_b(\R^d)$ with $\|f\|_{\tt W^{-\dd,k}}\le 1$ and $t\in (0,\tau(\gg)\land \tau(\tt\gg))$,
\beg{align*}  & \big|\gg(P_t^{\tt\gg} f)-\tt\gg(P_t^{\tt\gg} f)\big|
 = \bigg|\int_{\R^d\times\R^d} \big(P_t^{\tt\gg} f(x)- P_t^{\tt\gg} f(y)\big)\pi(\d x,\d y)\bigg|\\
&\le k_1  \bigg|\int_{\R^d\times\R^d} |x-y|\big(\scr M|\nn P_t^{\tt\gg} f|(x)+\scr M|\nn P_t^{\tt\gg} f|(y)+ {\|P_t^{\tt\gg} f\|_\infty}\big)\pi(\d x,\d y)\bigg|\\
&\le  k_1 {2^{\frac{1}{q}}}   \W_q(\gg,\tt\gg)  \Big[(\gg+\tt\gg)\big((\scr M|\nn P_t^{\tt\gg }f|)^{\ff{q}{q-1}}\big)\Big]^{\ff {q-1}q}+ {k_1\W_1(\gg,\tt\gg)\|P_t^{\tt\gg} \|_{\tilde{W}^{-\delta,k}\to \tilde{W}^{0,\infty}}}\\
&\le  k_1 {2^{\frac{1}{q}}}    \W_q(\gg,\tt\gg)   (\|\gg\|_{ {p_0}*}+\|\tt\gg\|_{ {p_0}*})^{\ff{q-1}q} \big\|\scr M|\nn P_t^{\tt\gg} f|\big\|_{\tt L^{k_0}}+ {k_1\W_1(\gg,\tt\gg)\|P_t^{\tt\gg} \|_{\tilde{W}^{-\delta,k}\to \tilde{W}^{0,\infty}}}\\
&\le k_2 \W_q(\gg,\tt\gg)   (\|\gg\|_{\vv,p*}+\|\tt\gg\|_{\vv,p*})^{\ff{q-1}q}  \big\| \nn P_t^{\tt\gg} \big\|_{\tt W^{-\dd,k}\to \tt W^{0,k_0}}+ {k_1\W_1(\gg,\tt\gg)\|P_t^{\tt\gg} \|_{\tilde{W}^{-\delta,k}\to \tilde{W}^{0,\infty}}}\\
&\le c_1(t) \W_q(\gg,\tt\gg)   (\|\gg\|_{\vv,p*}+\|\tt\gg\|_{\vv,p*})^{\ff{q-1}q}   \exp\big[t\bb_t k_t(\tt\gg)^\theta\big] t^{-[(\ff{1+\xi(q)}2)\vee(\frac{\delta}{2}+\frac{d}{2k})]}.\end{align*}
Therefore,
 \beq\label{X3} \beg{split} &\big\|P_t^{\tt \gg*}\tt\gg- P_t^{\tt\gg *}\gg\big\|_{\delta,k*}  \le D(t,\gg,\tt\gg)  t^{-[(\ff{1+\xi(q)}2)\vee(\frac{\delta}{2}+\frac{d}{2k})]},\\
&D(t,\gg,\tt\gg):=   c_1(t)   (\|\gg\|_{\vv,p*}+\|\tt\gg\|_{\vv,p*})^{\ff{q-1}q} \e^{t\bb_t  k_t(\tt\gamma)^{\theta}  }\W_q(\gg,\tt\gg),\\
&\qquad   \gg,\tt\gg\in \scr P_{\vv,p*},\ t\in (0,\tau(\gg)\land \tau(\tt\gg)).\end{split}\end{equation}

(b)  By Duhamel's formula  \cite[Proposition 5.5 (2)]{HRW25}, we have
$$P_t^{\gg}f-P_t^{\tt\gg}f= \int_0^t P_s^\gg \big\<b_s(\cdot,P_s^*\gg)- b_s(\cdot,P_s^*\tt\gg),\ \nn P_{s,t}^{\tt\gg}f\big\>\d s,\ \ f\in C_0^\infty(\R^d),$$
and   {\bf(A)} implies
$$ |b_t(x,P_t^*\gg)-b_t(x,P_t^*\tt\gg)|\le K_tt^\kk\|P_t^*\gg-P_t^*\tt\gg\|_{\dd,k*},\ \ t\in [0, \tau(\gg)\land  \tau(\tt\gg)).$$
 Then
\beq\label{X*2} \beg{split}
& \big\|P_t^{\gg *}\gg
 -  P_t^{\tt\gg *}\gg \big\|_{\dd,k*}=\sup_{f\in \B_b(\R^d), \|f\|_{\tt W^{-\dd,k}}\le 1}
 \big| \gg\big(P_t^{\gg}f-P_t^{\tt\gg}f\big)\big|\\
&\le K_t \|\gg\|_{\vv,p*} \int_0^t s^\kk  \big\|P_s^*\gg-P_s^*\tt\gg\big\|_{\dd,k*} \|P_s^\gg\|_{\tt W^{-\dd',k'}\to\tt W^{-\vv,p}}
\|\nn P_{s,t}^{\tt\gg}\|_{\tt W^{-\dd,k}\to\tt W^{-\dd',k'}}\d s,\\
&\qquad  \ \  t\in (0,\tau(\gg)\land \tau(\tt\gg)),\end{split}\end{equation}
where $(\dd',k')$ is in \eqref{TA} for $r\in [0,1]$ to be determined.

By \eqref{TJ}, we have either $\eta<\ff 1 2 +\kk$ and $\kk\ge \ff 1 2$, or  {$\eta\in [0,1)$} and $\kk<\ff 1 2$.
Let
\begin{equation*} r=\beg{cases} 1- \kk, &\text{if}\ \kk<\ff 1 2, \eta\in [0,1),\\
\ff{2}{3+2\kk}, &\text{if} \ \eta<\ff 1 2 +\kk, \kk\ge \ff 1 2.\end{cases}\end{equation*}
In each  case we have
\beq\label{YY0} r\eta <1,\ \ \ (2-r)\eta<2+2\kk,\ \ \ (1-r)\eta\le 2\kk,\end{equation}
where $(2-r)\eta<2+2\kk$ implies
$$\dd'-\vv+\ff{d(p-k')}{pk'}=(1-r) \eta<1\lor (2+2\kk-\eta),$$ so that Lemma \ref{L3'}(1) for $i=0$ gives
$$\|P_t^{\gg}\|_{\tt W^{-\dd',k'}\to\tt W^{-\vv,p}}\le \bb_t\e^{t\bb_t k_t(\gg)^\theta} t^{-\ff {1-r}2\eta},$$
and    { by Lemma \ref{L3'}(1) for $i=1$,}
$$\dd-\dd'+\ff{d(k'-k)}{k'k}=r\eta <1$$
implies
$$ \|\nn P_{s,t}^{\tt\gg}\|_{\tt W^{-\dd,k}\to\tt W^{-\dd',k'}}\le \bb_t\e^{t\bb_t k_t(\tt\gg)^\theta} (t-s)^{-\ff 1 2-\ff {r}2\eta}.$$
 Combining this with \eqref{X*2},  we find $c_2(t)\in (0,\infty)$ such that
  \beg{align*}& {\big\| P_t^{\gg *}\gg-P_t^{\tilde{\gg}*}\gg\big\|_{\dd,k*}} \le \tt D(t,\gg,\tt\gg)  \int_0^t s^{\kk-\ff{1-r}2\eta } \big\|P_s^*\gg-P_s^*\tt\gg\big\|_{\dd,k*} (t-s)^{-\ff 1 2-\ff{r}2\eta} \d s,\\
& \tt D(t,\gg,\tt\gg):=  { c_2(t) \|\gamma\|_{\vv,p*}   \e^{t\bb_t [k_t(\gamma)+k_t(\tt\gg)]^{\theta}},   \   \  t\in [0, \tau(\gg)\land  \tau(\tt\gg))}.\end{align*}
This together with \eqref{EQ} and   \eqref{X3} yields
\beg{align*}  &\big\|P_t^*\gg-P_t^*\tt\gg\big\|_{\dd,k*}\le  {\big\|P_t^*\gg- P_t^{\tilde{\gg} *}\gg\big\|_{\dd,k*}
+\big\|  P_t^{\tilde{\gg} *}\gg-P_t^*\tt\gg\big\|_{\dd,k*}}\\
&\le  D(t,\gg,\tt\gg)    t^{-[(\ff{1+\xi(q)}2)\vee(\frac{\delta}{2}+\frac{d}{2k})]}    +  \tt D(t,\gg,\tt\gg)   \int_0^t s^{\kk-\ff{1-r}2\eta} \big\|P_s^*\gg-P_s^*\tt\gg\big\|_{\dd,k*} (t-s)^{-\ff 1 2-\frac{r}2\eta} \d s,\\
&\qquad \gg,\tt\gg\in \scr P_{\vv,p*}, \  t\in [0, \tau(\gg)\land  \tau(\tt\gg)).\end{align*}
Since $(P_t^*\gg)_{t\in [0,T]}, (P_t^*\tt\gg)_{t\in [0,T]}\in \C_{\vv,p;\dd,k}^T$ for $T<\tau(\gg)\land\tau(\tt\gg)$, and
$1+\xi(q)\ge \eta$ due to \eqref{XQ},  for any constant $\ll\in (0,\infty)$ and  $t\in [0, \tau(\gg)\land  \tau(\tt\gg)),$   we have
\begin{equation*}I_t:=\sup_{s\in (0, t]} \e^{-\ll s}  s^{[(\ff{1+\xi(q)}2)\vee(\frac{\delta}{2}+\frac{d}{2k})]}\big\|P_s^*\gg-P_s^*\tt\gg\big\|_{\dd,k*}  <\infty\end{equation*}
and
 \begin{align*} I_t \le &\,D(t,\gg,\tt\gg) + \tt D(t,\gg,\tt\gg)    \sup_{s\in (0, t]}s^{[(\ff{1+\xi(q)}2)\vee(\frac{\delta}{2}+\frac{d}{2k})] + ([(\ff{1+\xi(q)}2)\vee(\frac{\delta}{2}+\frac{d}{2k})] +\ff{1-r}2\eta-\kk)^-}\\
 &\qquad\qquad\qquad\quad\quad\times
\int_0^s u^{  -([(\ff{1+\xi(q)}2)\vee(\frac{\delta}{2}+\frac{d}{2k})] +\ff{1-r}2\eta-\kk)^+} \e^{-\lambda(s-u)}(s-u)^{-\ff 1 2-\ff r 2 \eta}\d u.
\end{align*}
By \eqref{YY0}, \eqref{XQ} and $\dd+\frac{d}{k}< 2+(2\kk-\eta)^+$, we have
$$(\ff{1+\xi(q)}2)\vee(\frac{\delta}{2}+\frac{d}{2k})+ \ff{1-r}2\eta-\kk<1,\ \ \  \ff 1 2+\ff r 2 \eta<1,$$
we may apply  \eqref{LN0} to $\aa_1 =([(\ff{1+\xi(q)}2)\vee(\frac{\delta}{2}+\frac{d}{2k})] +\ff{1-r}2\eta-\kk)^+,\aa_2= \ff 1 2 +\ff r 2 \eta$ and $\aa= \ff 1 \theta,$ such that
$$I_t \le D(t,\gg,\tt\gg) + c\tt D(t,\gg,\tt\gg)   \ll^{-\ff 1 \theta}$$
holds for some constant $c\in (0,\infty)$.
Taking
$$\ll= \big(2c\tt D(t,\gg,\tt\gg) \big)^{\theta}$$
we obtain
$$I_t\le 2 D(t,\gg,\tt\gg),$$
which implies \eqref{ES5} for some increasing $\bb: (0,\infty)\to (0,\infty).$

Finally, when $\vv=0,p=\infty$, $k_t(\gg)$ defined in \eqref{1ga} is bounded above by some constant $c(t)\in (0,\infty)$ uniformly in $\gg\in \scr P_{\vv,p*}=\scr P$. Moreover,  in this case
\eqref{TJ} implies that $q=1$ satisfies \eqref{QY}.
So, \eqref{ES5} implies  \eqref{ES5'}.
\end{proof}

\beg{proof}[\textbf{Proof of Theorem \ref{T02}(2)}]  The proof is similar to that of \cite[Theorem 2.3(2)]{HRW25}, where $\vv=\dd=0.$ For completeness we figure out it in the present situation.

(a)  For $t\in (0, \tau(\gg)\land\tau(\tt\gg))$, we consider the SDEs
  \beq\label{LG} \beg{split} &\d X_s= b_s(X_s,P_s^*\gg)\d s + \d W_s,\\
  &\d Y_s= b_s(Y_s,P_s^*\tt\gg)\d s + \d W_s,\ \ \ s\in [0,t],\end{split}\end{equation}
where $X_0,Y_0$ are $\F_0$-measurable   such that
  \begin{equation*} \L_{X_0}=\gg,\ \ \  \L_{Y_0}=\tt\gg,\ \ \ \E|X_0-Y_0|^2= \W_2(\gg,\tt\gg)^2.\end{equation*}
Then
\beq\label{GT}   P_t^*\gg=\L_{X_t},\ \ \  \ P_t^*\tt\gg=\L_{Y_t}.\end{equation}
To estimate $\Ent(P_t^*\gg|P_t^*\tt\gg)$,  we
  apply the bi-coupling argument developed in \cite{23RW}.

For fixed $\theta'\in (\theta,\infty)$, let
\begin{align}\label{key}
t_0=\frac{t}{2}\wedge  s_t(\theta',\gamma)\wedge s_t(\theta',\tilde{\gamma}).
\end{align}
According to the bi-coupling method developed in \cite{23RW}, the following SDE  will be coupled with those two SDEs in \eqref{LG} respectively:
{ \begin{equation*}   \d Z_s =  \Big(1_{[0,t_0]}(s) b_s(Z_s, P_s^\ast\tilde{\gamma})+ 1_{(t_0,t]}(s) b_s(Z_s,P_s^\ast\gamma)\Big)\d s+ \d W_s,\ \ Z_0=Y_0,s\in [0,t].\end{equation*}}
 By \eqref{GT} and     \cite[Lemma 2.1]{23RW},   we have
\beq\label{RW} \beg{split}
&\Ent(P_t^\ast\gamma| P_t^\ast\tilde{\gamma})=\Ent(\L_{X_t}|\L_{Y_t})\\
&\le 2 \Ent(\L_{X_t}|\L_{Z_t}) +  \log \int_{\R^d} \Big(\ff{\d\L_{Z_t}}{\d\L_{Y_t}}\Big)^{2}\d \L_{Y_t}
  =:2 I_1+  I_2.\end{split}\end{equation}
 Below we estimate $I_1$ and $I_2$  respectively.

(b) To estimate $I_1$, we   first establish the   log-Harnack inequality  for $P_{t}^\gg:$ for any  {$\theta'\in (\theta,\infty)$,}  there exists   $c_0 (t)  \in (0,\infty)$ increasingly in $t$  such that
\beq\label{cty}\beg{split}
&P_{s,t}^\gamma \log f(x) \leq \log P_{s,t}^\gamma f(y)+\frac{c_0(t) |x-y|^2}{s_t(\theta',\gg)\wedge (t-s)},\\
&\qquad  \ x,y\in\R^d,   0\leq s< t<\tau(\gg),\    \gg\in \scr P_{\vv,p*},\  f\in \B_b^+(\R^d)
\end{split}\end{equation}
 for $s_t(\theta',\gg)$ defined in \eqref{ST}. By {\bf (A)}, we have
 $$\|b_t(\cdot,P_t^\ast\gamma)\|_\infty\le K_t k_t(\gg)t^{\kk-\ff\eta 2},\ \ \ t\in [0,\tau(\gg)).$$
 Since \eqref{TJ} implies \eqref{TJ'} so that $(\eta-2\kk)^+<1$. For any constant
 $$\theta'> \theta:=\ff 2 {1-(\eta-2\kk)^+},$$ we have
 $$q':=  {\Big(\ff{(\eta-2\kk)^+}2 + \ff 1 {\theta'}\Big)^{-1} }\in \Big(2, \ff 2 {(\eta-2\kk)^+}\Big),$$ and for some $c_1(t)\in(0,\infty)$ increasing in $t$,
\beq\label{GM}\|b_\cdot(\cdot,\gg_\cdot)\|_{\tt L_{q'}^\infty(s,t)}\le c_1(t) k_t(\gg)(t-s)^{\ff 1 {\theta'}},\ \ 0\le s<t<\tau(\gg).\end{equation}
  By \eqref{ST},  we find a constant $k_1\in (0,\infty)$ such that
 $$k_t(\gg)(t-s)^{\ff 1{\theta'}}\le k_1,\ \ 0<t-s\le  s_t(\theta', \gg).$$
 {This together with \eqref{GM} implies that}
  $$\|b_\cdot(\cdot,\gg_\cdot)\|_{\tt L_{q'}^\infty(s,t)} \le  {k_1c_1(t)},\ \ 0<t-s\le s_t(\theta',\gg), \ t\in (0,\tau(\gg)).$$
So, by \cite[Proposition 5.2(4)]{HRW25},  we derive \eqref{cty}  for $0<t-s\le s_t(\theta',\gg), \ t\in (0,\tau(\gg))$ and some    $c_0 (t)  \in (0,\infty)$ increasingly in $t$.
When   $t-s>s_t(\theta',\gg),\ t\in (0,\tau(\gg)),$  by the semigroup property and Jensen's inequality,  we deduce
  \beg{align*}&
 P_{s,t}^\gg \log f(x)=P_{s,s+s_t(\theta',\gg)}^\gg P_{s+s_t(\theta',\gg),t}^\gg \log f(x)\leq P_{s,s+s_t(\theta',\gg)}^\gg \log P_{s+s_t(\theta',\gg),t}^\gg  f(x)\\
 &\leq \log P_{s,s+s_t(\theta',\gg)}^\gg P_{s+s_t(\theta',\gg),t}^\gg  f(y)+\ff{ {c_0(t)|x-y|^2}}{s_t(\theta',\gg)}\\
  &= \log P_{s,t}^\gg f(y)+ \ff{ c_0(t)|x-y|^2}{s_t(\theta',\gg)}, \ \ x,y\in\R^d,\ f\in \B_b^+(\R^d).\end{align*}
 So, \eqref{cty} also  holds for  $t-s>s_t(\theta',\gg)$.

 Next, by the Markov property, we have
\begin{equation*}\E[f(X_t)]= \E[ (P_{t_0,t}^\gg f)(X_{t_0})],\ \ \
\E[f(Z_t)] = \E[ (P_{t_0,t}^\gg f)(Z_{t_0})].\end{equation*}
This together with \eqref{cty} for $s=t_0$ and Jensen's inequality implies
\beq\label{OG3}\E[\log f(X_t)]\le \log \E[f(Z_t)]+ \ff {2c_0(t)}{  s_t(\theta',\gg)}\E[|X_{t_0}- Z_{t_0}|^2],\ f\in \B_b^+(\R^d).\end{equation}
Moreover,  by  {\bf (A)},  {\eqref{ES5}} and \eqref{key},  we find $c_2(t) \in (0,\infty)$ increasing in $t$ such that
 \beg{align*} &\int_0^{t_0} \|b_s(\cdot,P_s^\ast\gamma)-b_s(\cdot,P_s^\ast\tilde{\gamma})\|_\infty\d s  \le K_t \int_0^{t_0} \|P_s^\ast\gamma-P_s^\ast\tilde{\gamma}\|_{\dd,k} s^\kk\d s \\
 &\le K_t \big(\|\gg\|_{\vv,p*}+\|\tt\gg\|_{\vv,p*}\big)^{\ff{q-1}q} K_{t,\bb}^{(\theta)}(\gg,\tt\gg)   \W_q(\gg,\tt\gg)\int_0^{t_0}  s^{\kk- [(\ff{1+\xi(q)}2)\vee(\frac{\delta}{2}+\frac{d}{2k})]}\d s\\
 &\le c_2(t) \big(\|\gg\|_{\vv,p*}+\|\tt\gg\|_{\vv,p*}\big)^{\ff{q-1}q}    \W_q(\gg,\tt\gg) t_0^{ {(1+\kappa)- [(\ff{1+\xi(q)}2)\vee(\frac{\delta}{2}+\frac{d}{2k})]}}. \end{align*}
 Combining this with \cite[Proposition 5.5(1)]{HRW25},  we find $c_2(t) \in (0,\infty)$ increasing in $t$ such that
 $$\E[|X_{t_0}- Z_{t_0}|^2] \le c_2(t) \big(\|\gg\|_{\vv,p*}+\|\tt\gg\|_{\vv,p*}\big)^{\ff{2(q-1)}q}  \big(\W_q(\gg,\tt\gg)^2  {t_0}^{2(1+\kappa)-[(1+\xi(q))\vee(\delta+\frac{d}{k})]}+ \W_2(\gg,\tt\gg)^2\big).$$
 This together with  the formula
$$\Ent(\mu|\nu)=\sup_{f\in \B_b^+(\R^d)} \big(\mu(\log f)-\log \nu(f)\big),\ \ \mu,\nu\in \scr P $$
 and \eqref{OG3} yields that for some  constant $c_3(t) \in (0,\infty)$ increasing in $t $ such that
 \beq\label{I11}\begin{split}
 &I_1:=\Ent(\L_{X_t}|\L_{Z_t})=\sup_{f\in \B_b^+(\R^d)} \big(\E[\log   f(X_t)]-\log \E[ f(Z_t)]\big),  \\
 & \leq   c_3(t)   (\|\gg\|_{\vv,p*}\lor\|\tt\gg\|_{\vv,p*})^{\ff{2(q-1)}q}\bigg( \ff{\W_2(\gg,\tt\gg)^2}{ s_t(\theta',\gamma)}+ \ff{\W_q(\gg,\tt\gg)^2}{s_t(\theta',\gg)^{([(1+\xi(q))\vee(\delta+\frac{d}{k})]-(2\kappa+1) )^+}}\bigg),
 \end{split}\end{equation}
where in the last step we have used $t_0\le   s_t(\theta',\gg)$.

(c) Estimate $I_2$.
  By {\bf(A)} and  \eqref{ES5},  we find a constant $K(t)\in (0,\infty)$  increasing in $t$   such  that
\begin{equation*}\xi_s := \big[b_s(Y_s, P_s^\ast\gamma)- b_s(Y_s, P_s^\ast\tilde{\gamma})\big]\end{equation*} satisfies
\beq\label{XX}\beg{split}& |\xi_s |^2 \le K_t^2 C(t,\gg,\tt\gg)s^{2\kk-[(1+\xi(q))\vee(\delta+\frac{d}{k})]},\ \ \ s\in (0,t],\\
&C(t,\gg,\tt\gg):=   (\|\gg\|_{\vv,p*}+\|\tt\gg\|_{\vv,p*})^{\ff{2(q-1)}q}  K_{t,\beta}^{(\theta)}(\gg,\tt\gg)^2 \W_q(\gg,\tt\gg)^2.\end{split}\end{equation}
 Then
 $$ R_s:=\e^{\int_{t_0}^s \<\xi_r,\d W_r\>-\ff 1 2 \int_{t_0}^s|\xi_r|^2\d r},\ \ s\in [t_0,t]$$
is a martingale, and by Girsanov's theorem,
$$\ff{\d\L_{Z_t}}{\d\L_{Y_t}}(Y_{t})= \E(R_t|Y_{t}).$$
Combining this with Jensen's inequality and \eqref{XX},
we find   $c_4(t)\in (0,\infty)$  increasing in $t$  such that
   \beg{align*}I_2&:= \log \E\bigg[\Big(\ff{\d\L_{Z_t}}{\d\L_{Y_t}}(Y_{t})\Big)^{2}\bigg]
\leq \log\E\Big[R_t^{2}\Big]\\
&\le \log\E\bigg[\e^{2\int_{t_0}^t \<\xi_s,\d W_s\>- 2 \int_{t_0}^t |\xi_s|^2\d s+  C(t,\gg,\tt\gg) \int_{t_0}^t s^{2\kk- 1-\xi(q)}\d s }\bigg]\\
&=  C(t,\gg,\tt\gg) \int_{t_0}^t s^{2\kappa-[(1+\xi(q))\vee(\delta+\frac{d}{k})]}\d s
 \le  c_3(t) C(t,\gg,\tt\gg) t_0^{-([(1+\xi(q))\vee(\delta+\frac{d}{k})]-(2\kk+1))^+}.
\end{align*}
By combining this with   \eqref{RW} and \eqref{I11}, we obtain \eqref{ES6} for some    $\bb: (0,\infty)\to (0,\infty).$

(c) If $\vv=0,p=\infty$, we have   $\scr P_{\vv,p*}=\scr P$, $\tau(\gg)=\infty$ and $\|\gg\|_{\vv,p*}= 1$ for any $\gg\in \scr P$,
and we may take $q=1$ so that $\xi(q)=\eta=\delta+\frac{d}{k}$,
and \eqref{XQ} implies  $$[(1+\xi(q))\vee(\delta+\frac{d}{k})]-(2\kk+1)  <1.$$ Hence \eqref{ES6} implies \eqref{ES7}
for some increasing $\bb: (0,\infty)\to (0,\infty).$
\end{proof}



\end{document}